\documentclass[aip,preprint]{revtex4-2}

\usepackage{siunitx}
\usepackage{amsfonts,amssymb,amsthm,amsmath}
\usepackage{mathrsfs}
\usepackage{mathtools}
\usepackage{physics}
\usepackage{mathrsfs}
\usepackage{marvosym}
\usepackage{mathtools}
\usepackage{subcaption}
\usepackage{caption}
\usepackage{lineno}

\newcommand{\rmd}{\mathrm{d}}

\draft 

\begin{document}


\title{Guidelines for data-driven approaches to study transitions in multiscale systems: the case of  Lyapunov vectors} 




\author{Akim Viennet}
\email[]{akim.viennet@ens.fr}
\affiliation{Department of Physics, Ecole Normale Superieure of Paris}

\author{Nikki Vercauteren}
\email[]{nikki.vercauteren@geo.uio.no}
\affiliation{Department of Geosciences, University of Oslo
}
\author{Maximilian Engel}
\email[]{maximilian.engel@fu-berlin.de}
\affiliation{Institute of Mathematics, Freie Universit{\"a}t Berlin
}
\author{Davide Faranda}
\email[]{davide.faranda@lsce.ipsl.fr}
\affiliation{Laboratoire des Sciences du Climat et de l'Environnement, UMR 8212 CEA-CNRS-UVSQ, Universit\'e Paris-Saclay, IPSL, 91191 Gif-sur-Yvette, France}
\altaffiliation{London Mathematical Laboratory, 8 Margravine Gardens London, W6 8RH, UK}
\altaffiliation{LMD/IPSL, Ecole Normale Superieure, PSL research University, Paris, France}

\date{\today}

\begin{abstract}
We study in detail the role of covariant Lyapunov vectors and their respective angles for detecting transitions between metastable states in dynamical systems, as recently discussed in several atmospheric science applications. The underlying models are built from data by the dynamical clustering method, called FEM-BV-VAR, and the Lyapunov vectors are approximated based on these models. We test this data-based numerical approach at the hand of three well-understood example systems with increasing dynamical complexity, identifying crucial properties that allow for a successful application of the method: in particular, it turns out that the method requires a clear multiple time scale structure with fast transitions between slow subsystems which can be dynamically characterized by invariant neutral directions of the linear approximation model.
\end{abstract}


\maketitle 

\tableofcontents

\section{Introduction}\label{Sec:intro}

Dynamical systems theory deals with the prediction of trajectories of natural systems originating by one given or a set of initial conditions. This task is particularly challenging when the system is chaotic and even more when the system features several metastable states or multiscale features~\cite{katok1997introduction,manneville2010instabilities}.  Among all the possibilities, here we focus on the  stability properties of certain meta-stable states, that organize the phase space, and the estimation of the probability of switching from one state to another. Many mathematical tools have been developed to address those questions, one of them being the study of the so-called \emph{covariant Lyapunov vectors} (CLVs) (also known as Oseledets vectors), and the associated \emph{Lyapunov exponents} (LEs). These vectors give a basis on the tangent space at points of trajectories, providing directions of linear perturbation growth along the dynamics \cite{Ruelle_1979, Ginelli_2007, Wolfe_2007}. They can be seen as a generalization of the linear stability theory for fixed points and of Floquet's theory for limit cycles, since Lyapunov vectors and exponents can be computed along any trajectory of a smooth dynamical system. The CLVs give the directions of growth or decay of a perturbation, and the LEs give the associated rate of asymptotic growth or decay. An increase of one of the unstable LEs has been associated to a higher instability for various theoretical and physical systems \cite{Gilmore_2019, Nazarimehr_2017}. 

In many cases, transient (chaotic) behavior cannot be detected by asymptotic LEs which average out transient dynamics via ergodic limits; hence, \emph{finite-time Lyapunov exponents} (FTLEs) are often more suitable to capture the degree of uncertainty at different points of trajectories and their small neighbourhoods.

The directions of unstable CLVs indicate the directions towards which an error will grow with the rates given by the associated (FT)LEs. For example, this tool can be used in ensemble weather forecasting to identify how to enforce initial perturbations to optimally span the space of possible realizations of the weather \cite{Toth_1993}.
Another quantity of interest is the angle between the flow direction and the most unstable CLV. An alignment of those vectors has been used as a predictor for transitions, tipping points or catastrophes (extreme events) in several systems \cite{Sharafi_2017, Beims_2016}. In particular, this criterion has been proven to be an important early-warning sign for abrupt transitions in the Peña and Kalnay climate toy-model\cite{Quinn_2020}.
Finally, Quinn et al.\cite{Quinn_2021} suggested that the projection of the most unstable CLV just before a transition between two states could inform on the patterns that triggered the instability and then the transition.

Summarizing, the computation of CLVs and associated LEs is of high interest for the analysis of dynamical systems. Recent progress was made to compute them numerically, due to various algorithms by Ginelli et al. \cite{Ginelli_2007}, Wolfe and Samelson \cite{Wolfe_2007}, and Froyland et al. \cite{Froyland_2013}. However, all those algorithms rely on the knowledge of an analytic expression of the model, in order to differentiate the flow and compute the linear cocycles (see Section~\ref{CLVs algo} for an introduction to those methods). This suggests that it is rather difficult to use such methods based on observations for which the underlying model is unknown or only partially known, such as reanalysis atmospheric data.

Yet, Quinn et al. \cite{Quinn_2021} recently introduced a method to compute CLVs directly from data, by introducing a model-based clustering step before estimating the CLVs. A model is fitted to the observations via a dynamical, or model-based clustering method initially introduced by Horenko \cite{Horenko_2010, Metzner_2012}: the FEM-BV-VAR (Finite Element clustering with bounded variation (FEM BV) vector autoregressive (VAR)) clustering approach. 
Differently from more classical, geometrical clustering methods, in this framework a state is not defined by a geometrical area in the phase space, but by an estimated auto-regressive dynamics (see Section~\ref{FEM} for details). 
The whole system is then switching between those dynamical models. 
This method is particularly adapted to the purpose since it provides not only a cluster affiliation sequence, but also a linear (auto-regressive) model for each of the states. One can then use these (approximated) models to compute an approximation of the CLVs and of the LEs for the dynamical system underlying the data, and thus get some insights on the stability of the states and of the stable and unstable directions. 
Quinn et al.\cite{Quinn_2021} used this approach to analyse the dynamics of atmospheric circulation patterns in the northern hemisphere. They investigated the dynamical stability properties of recurrent and persistent states of the atmospheric circulation patterns or regimes known as the North Atlantic Oscillation (NAO) and atmospheric blocking events. In particular, the CLVs were used to analyse the pressure distribution patterns related to transitions between the recurrent circulation regimes, leading to insightful observations since weather forecasting and climate models struggle to capture the onset and decay of blocking events.

These results led to the question whether the method is applicable for other systems of interest and to which extent it more generally captures relevant information on the dynamics. The aim of this work is to explore this question by testing the method in several systems for which some a priori knowledge of the CLVs and of the transitions between regimes is available: a fast-slow FitzHugh-Nagumo oscillator, a well-studied Von Kármán turbulent flow from a laboratory experiment, and a Lorenz 63 system, where the order of our presentation follows an increase of dynamical complexity.
Results on those different systems will show that the method provides several insights on the dynamics, quantifying the stability of different (meta-stable) states and thereby identifying transitions between them. 
This holds true in particular for the Von Kármán flow. 
However, we also demonstrate why such conclusions may be treated with caution, considering the strong dependence on the existence of an dynamically invariant normal tangent flow direction, a visible time scale separation and a large number of hyper-parameters.

In this paper, we will first present the details of the method that allows one to compute approximated CLVs from a data series. Then we will try to assess its validity, by applying it on a FitzHugh-Nagumo oscillator, experimental data from the Von Kármán flow and the Lorenz 63 model.
Finally, we will discuss the scope of the methods at the hand of these examples, illustrating its potential but also several caveats for its application.

\section{The Lyapunov vectors and their numerical computation}
\subsection{Mathematical background}

\begin{figure}[hb]
    \centering
    \includegraphics[width=\linewidth]{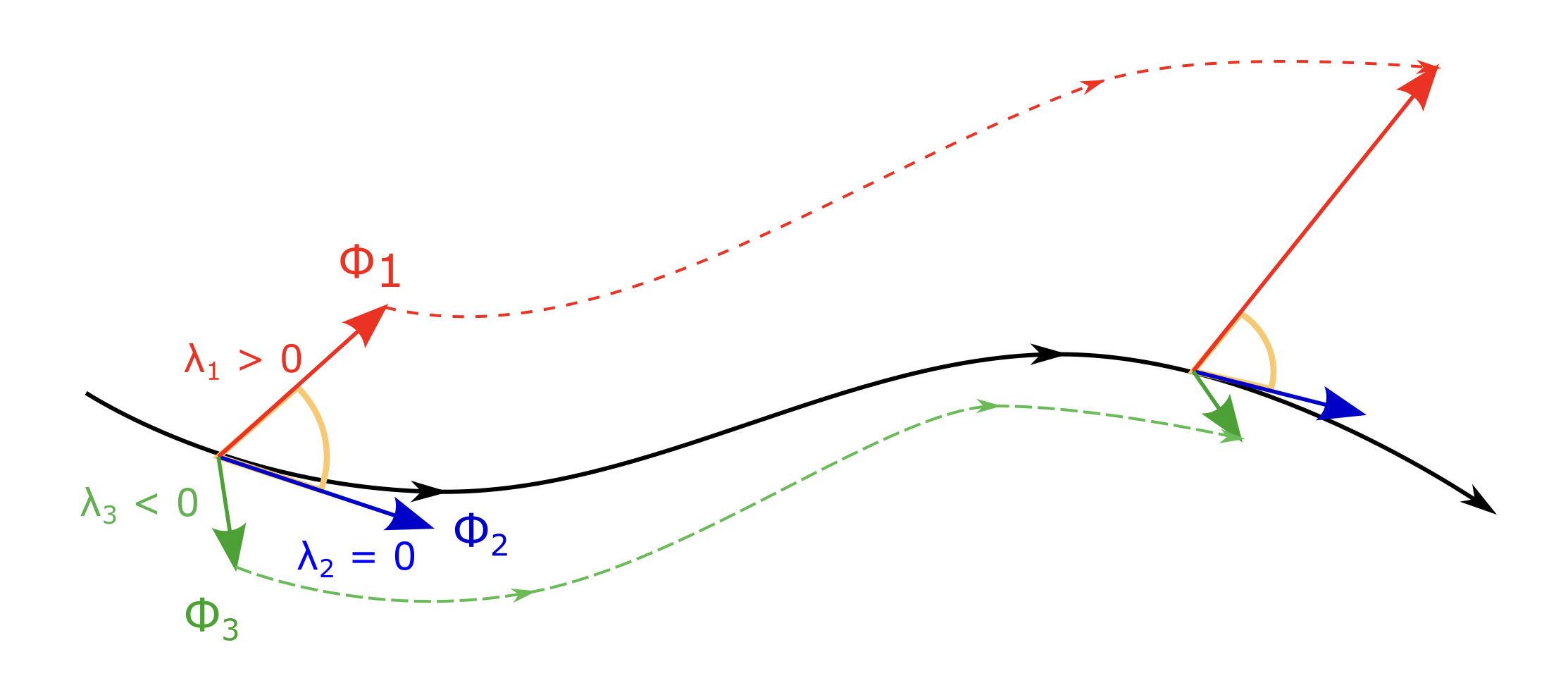}
    \caption{Contraction and expansion of CLVs along a trajectory with positive, negative and zero Lyapunov exponent $\lambda_i$, $i=1,2,3$. In this setting, we have at each point three CLVs ($\phi_1$, $\phi_2$ and $\phi_3$). The solid line represents the unperturbed trajectory, while the dotted lines represent the perturbed trajectories, along the stable (in green) and the unstable (in red) directions. The so-called alignment $\theta_{12}$ is given by the cosine of the orange angle.}
    \label{schema CLVs}
\end{figure}

Let us first introduce the notion of CLVs. They arise from a non-autonomous generalization of the linear stability analysis at fixed points and  Floquet theory at limit cycles to any point of the trajectory. For a dynamical system, the CLVs form a basis of the tangent space and give the directions of growth or decay of any perturbation around a background flow. The Lyapunov exponents (LEs) give the associated rate of growth or decay (see Fig.\ref{schema CLVs}). 
Assuming ergodicity of a dynamical system $\Phi_t(x_0)$, whose trajectories we will simply denote by $x(t)$, one observes that the LEs are global numbers 
that characterise the whole attractor, whereas the CLVs may depend on the particular points of the trajectory (but are still asymptotic objects). 

In more detail, the existence of Lyapunov exponents with corresponding directions on the tangent space is given by Oseledets' Multiplicative Ergodic Theorem (MET)\cite{MET}. Under a mild integrability condition with respect to an ergodic invariant measure,
this theorem gives us, in each point of the trajectory, the existence of a splitting of the tangent space into $p \leq d$ subspaces 
$$\mathbb{R}^{d}=Y_{1}(x(t)) \oplus \cdots \oplus Y_{p}(x(t)),$$
such that for all $ v \in Y_{i}(x(t))$,
\begin{equation}
\label{lim MET}
\lim _{\tau \rightarrow \infty} \frac{1}{\tau} \log \|\mathcal{F}\left(t, t+\tau\right) \cdot v \| = \lambda_{i},
\end{equation}
where $\mathcal{F}$ denotes the linear propagator for the tangent flow, i.e.
$$
\boldsymbol{v}\left(t_{2}\right)=\mathcal{F}\left(t_{1}, t_{2}\right) \boldsymbol{v}\left(t_{1}\right)
$$
and
$\lambda_{1}> \lambda_{2} > \cdots > \lambda_{p}$ are the distinct LEs with multiplicities $m_i \geq 1$, $i=1, \dots, p$. 
The CLVs $v_i^j(t)$, $j=1,\dots, m_i$, are then representative vectors from the Oseledets subspaces $Y_{i}(x(t))$, which are unique up to scalar factors if $m_i=1$ and chosen as a set of $m_i$ linearly indepedent vectors in $Y_{i}(x(t))$ otherwise. Let us order them as $\phi_k$, $k=1, \dots, d$, where $v_1^1= \phi_1, \dots, v_1^{m_1} = \phi_{m_1}$, and so on (see also Figure~\ref{schema CLVs} where all $m_i=1$ as will be the case in our examples).

As described in Section~\ref{Sec:intro}, LEs and CLVs give important information on the stability properties of the dynamics, and have been used to predict transitions and extreme events. One key quantity is the angle between the neutral CLV (the CLV associated with a zero LE, which is always tangent to the flow direction) and the most unstable CLV (the one associated with the largest positive LE), given that they both exist. Let us call $\theta_{ij}$ the cosine of the angle between the CLVs $\phi_i$ and $\phi_j$:

\begin{equation}
    \theta_{ij}(t) = \frac{\abs{\phi_i(t) \cdot \phi_j(t)}}{\lVert \phi_i(t)\lVert \dot \lVert \phi_j(t) \lVert}
    \label{def theta}
\end{equation}
%
Many studies suggest that, for $\phi_i$ representing the most unstable direction and $\phi_j$ a neutral direction, this angle is related to the probability of transitions between characteristic states: the more these two vectors align, the higher such a switching probability is expected to be. Sharafi et al.\cite{Sharafi_2017} have applied this criterion to various fast-slow systems, whereas Beims et al.\cite{Beims_2016}  have used it to predict extreme events in a Rössler oscillator. In the following, we will call “alignment of CLVs” the absolute value of the cosine of the angle between a most unstable CLV and a neutral one. 
In cases without a neutral direction, one may take the CLV associated with the Lyapunov exponent closest to $0$ and consider this direction as a \emph{near-neutral} one (see also Section~\ref{sec:Lorenz}).

Note that the CLVs and associated Lyapunov exponents are asymptotic objects whereas the transitions we are interested in happen on finite time scales and the analyzed time series are also naturally finite. Hence, in reality one analyzes \emph{Finite Time Lyapunov Exponents} (FTLEs) which are defined analogously for a given finite $\tau$ in Eq.\eqref{lim MET}, depending on space and time. FTLEs associated to CLVs (or their finite time approximations one may also regard to as Finite Time Lyapunov vectors (FTLVs)) may change their signs depending on $\tau$. In some cases, there can exist strictly positive FTLEs even though the trajectory is asymptotically stable (i.e.~all LEs are negative). This is typical for globally asymptotically stable systems with transient chaos.
In particular, an asymptotically stable (or unstable) CLV might be referred to as an unstable (or stable) CLV on certain finite time scales. Hence, we will call CLVs stable or unstable in our numerical studies based on the local stability within the investigated finite time scales. 
The FitzHugh-Nagumo oscillator discussed below exemplifies this: while the trajectories asymptotically approach a stable periodic orbit (one negative, one neutral LE), the CLVs and associated FTLEs can detect the local instability along the fast subsystem (see Section \ref{sec:FHN} for further explanations). 
Generally speaking, FTLEs can be used as a measure for the predictability of the local dynamics: the higher the largest FTLE, the lower the predictability on the respective time scale (see, for example, Deremble et al.~\cite{Deremble_2009} for the classical Lorenz 63 attractor and a one-layer quasi-geostrophic atmospheric model).  As suggested by Quinn et al.\cite{Quinn_2021} (and before by Deremble et al.~\cite{Deremble_2009}), the time length $\tau$ acts as a scale filter for the dynamics: with small $\tau$ the computed FTLEs and the related CLVs (or FTLVs) give insights on the short scale processes, whereas with larger $\tau$ we get closer and closer to asymptotic properties. 

\subsection{Direct computation of the CLVs}\label{CLVs algo}

There exist several algorithms to numerically compute the CLVs. One of the most famous methods was developed in 2007 by Ginelli et al.\cite{Ginelli_2007}. However, here we will use a modified approach introduced in 2013 by Froyland et al.\cite{Froyland_2013} (algorithm 2.2 in this reference). This choice is motivated empirically by a faster convergence and by more consistent results in the considered setting, when compared to results obtained with Ginelli's algorithm. 

The  \emph{Froyland algorithm} relies on a singular value decomposition of the forward cocyles starting at past fibers, then propagating the obtained orthogonal directions into covariant ones. Thus, computing the CLVs at a given point on the trajectory requires a pullback procedure from the past to the present (and beyond).
This involves a number of time steps $N$ for going to the past and a number of time steps $M$ corresponding with the time length $\tau$ in Eq.~\eqref{lim MET}. In this study, for simplicity we always take $M = N$ (as suggested in Froyland's article and validated empirically). In theory, increasing $N$ and $M$ improves the approximation. 
However, our results show that convergence may fail due to the accumulation of numerical errors. Therefore, $N$ and $M$  are key parameters that act as a scale filter, similarly to $\tau$ in the previous subsection.
Another internal parameter to be adapted is given by the number of correction steps $n$ for obtaining the covariant out of singular directions; for details see algorithm 2.2 in Froyland et al.~\cite{Froyland_2013}. To sum up, this algorithm requires to set three parameters: 
$$
M, N \text{ and } n
$$
with, in this study, $N = M$.

Finally, let us emphasise that this algorithm requires an explicit expression for the linear propagator at each point. 
For continuous-time systems $\dot{x}=f(t, x)$, the linear propagator solves the variational linear differential equation with matrix generator $J(x, t) := \left(D_{x} f\right)(t, x)$, i.e.~the Jacobian of the vector field $f$. 
For discrete-time systems $x_{n+1} = g(x_n)$, the propagator is the product of the matrices $A_n= \left(D_{x} g\right)(x_n)$. 
Hence,  computing the quantities directly from data, for which the propagator is not known a priori, is out of reach. 
The aim of this article is to investigate the capabilities of the above-mentioned algorithm for computing approximate CLVs directly from observed time series, relying on a prior modelling step using a model-based clustering framework. We hence explore, based on systems of different complexity, the conditions under which the method first introduced in \cite{Quinn_2021} provides reliable results. 

\section{Dynamical clustering method}

\subsection{FEM-BV-VAR approach}\label{FEM}

In the literature, various approaches address the problem of identifying persistent states based on data. They can be roughly classified as either \textit{non-dynamical} or \textit{dynamical} methods. The class of non-dynamical methods only exploits geometrical properties of the data for clustering, regardless of their temporal occurrence. The most used non-dynamical approach is the k-means method, which clusters data points according to their minimal distance to geometrical centroids of point clouds \cite{hartigan_algorithm_1979}. Dynamical methods additionally take into account the temporal changes of data, based on latent variables models such as hidden Markov models \cite{rabiner_tutorial_1989}. This work considers a dynamical clustering method in which the existence of multiple states is presumed, each having time-independent properties. Those states are presumed to have a certain degree of persistence, and the system transitions between them during its evolution. A simplified description of the dynamics is then given in terms of a set of locally stationary linear vector autoregressive models (the cluster states). This method is coined as FEM-BV-VAR approach (Finite Element clustering with bounded variation (FEM BV) Vector autoregressive (VAR)) \cite{Horenko_2010, Metzner_2012}. Due to its proven utility in modeling transitional behavior between persistent meta-stable states directly from data, FEM-BV-VAR has recently become popular to study dynamical aspects of the atmosphere, ocean, and climate systems; studies have tackled small-scale processes in the atmospheric boundary layer \cite{vercauteren_clustering_2015, boyko_multiscale_2021}, as well as large-scale atmospheric and oceanic circulation \cite{FEM-OKane, okane_decadal_2013, Quinn_2021}. Importantly, the method does not rely on any underlying assumptions regarding the statistical stationarity of the data and, hence, is applicable to problems where trends are present.

In the FEM-BV-VAR approach, a cluster is defined as a subset of the observed time series of data whose evolution can be described approximately by a stationary linear vector autoregressive model. The full time series is modeled as a set of such stationary VAR models, with a switching process representing transitions between the cluster states. Since the states are assumed to have a certain degree of persistence, the dynamical evolution of the system is described by VAR models describing the fast-scale dynamics within a give state, while the slow evolution is described by the switching process. Hence, the dynamics is decomposed into two parts:
\begin{itemize}
    \item a locally stationary fast auto-regressive (VAR) process,
    \item a slow hidden process that makes the system switch between different forms of such auto-regressive processes (i.e. between the different states).
\end{itemize}

Within a given state, we assume the time evolution of the vector of observables  $\mathbf{x}_{t}$ to be governed by 
\begin{equation}\label{eq:VAR}
\mathbf{x}_{t}=\mu^{(i)}+\sum_{\tau=1}^{m} \mathbf{A}_{\tau}^{(i)} \mathbf{x}_{t-\tau}+\epsilon_{t}^{(i)}
\end{equation}
where $\mu^{(i)}$ is the mean of the $i$-th cluster, $\mathbf{A}_{\tau}^{(i)}$ are matrices, and $\epsilon_{t}^{(i)}$ is a white noise with a covariance matrix $\Sigma^{(i)}$ .
A state of the system (or cluster) $i$ is then characterized by its set of parameters
$$
\Theta_{i}=\left(\mu^{(i)}, \mathbf{A}_{1}^{(i)}, \ldots, \mathbf{A}_{m}^{(i)}, \Sigma^{(i)}\right).
$$
A set of $K$ such models is assumed, with different model coefficients in \eqref{eq:VAR}, leading to $K$ clusters. Determination of the optimal coefficients in \eqref{eq:VAR} is done via minimization based on the distance between the observations and the deterministic part of the model
\begin{equation}\label{eq:dist}
g \left(\mathbf{x}_{t},\mathbf{\theta}(t) \right) = \lVert \mathbf{x}_{t} - \mu^{(i)}(t) - \sum_{\tau=1}^{m} \mathbf{A}_{\tau}^{(i)} (t) \mathbf{x}_{t-\tau} \rVert,
\end{equation}
calculated for a fixed temporal realisation of parameters $\mathbf{\theta}(t)$.
The functional to minimize also includes a cluster affiliation term that determines the set of model parameters the data should be associated with and is then given as
\begin{equation}\label{eq:dist_func}
L\left(\mathbf{\Theta},\mathbf{\Gamma}(t)  \right) = \sum_{t=0}^{T} \sum_{i=1}^{K} \gamma_{i}(t) g \left(\mathbf{x}_{t},\Theta_{i} \right),
\end{equation}
where $\mathbf{\Theta}$ denotes the collection of all $\Theta_i$, i.e. $\mathbf{\Theta}=(\Theta_1, \cdots, \Theta_K)$ and $T$ the time length of the observed dynamics.
The functions $\mathbf{\Gamma}(t)=(\gamma_{1}(t), \cdots, \gamma_{K}(t))$ are the cluster affiliation functions whose values give the probability of the data at time $t$ to belong to cluster $i$ and should satisfy the following property at a given time $t$
\begin{equation}\label{eq:gamma}
\sum_{i=1}^{K} \gamma_{i} = 1, \quad \gamma_{i} \ge 0 \quad \forall i=1,\cdots, K
\end{equation}
The number $K$ and the memory depth $m$ are hyper-parameters that must be selected. The assumption of local stationarity of the statistical process is finally enforced by setting a persistence parameter $C$, which defines the maximum allowed number of transitions between a total of $K$ different statistical processes. This step regularises the minimization problem by introducing the additional constraint on the total variation norm of the sequence
\begin{equation}\label{eq:persis}
\sum_{t=0}^{T-1} |\gamma_{i}(t+1) - \gamma_{i}(t) | \le C, \quad \forall i=1,\cdots, K.
\end{equation}
This last hyper-parameter $C$ is also more conveniently defined via the average persistence $p$ as $C = \frac{T}{p} - 1$. The reader is referred to Horenko \cite{Horenko_2010} and references therein for further details about the method and the minimization process.  

This method makes it possible to detect dynamical patterns that would not be detected by a geometrical method such as the k-means: for instance, a change in frequency of the signal or some oscillations with multiple amplitudes. It also provides a local linear model for the data, on which the computation of the Covariant Lyapunov Vectors will be based. However, it is important to bear in mind that three hyper-parameters ($K, m, p$) have to be selected when fitting a model.


\subsection{Choosing the hyper-parameters}
\label{choice parameters FEM}
Statistical techniques based on information theory were developed to find the best hyper-parameters of the FEM-BV-VAR (namely the number of clusters $K$, the memory depth $m$ and the average persistence $p$)\cite{Horenko_2010, Metzner_2012}. Here, physical understanding of the systems is also used to choose $K$ and $m$, as will be detailed when presenting the results. The persistence $p$ is selected via the so called L-curve method: as shown by Horenko\cite{Horenko_2010},  the optimal value of $p$ can be determined as the edge point (or the point of maximal curvature) on a two- dimensional plot, where one plots the total distance between the model and the data against the value of $p$. In the application of the FEM-BV-VAR algorithm, the reconstructed signal has been found to diverge in some configurations; hence, we have checked the output of the algorithm manually and sometimes slightly modified $p$ around its optimal value if the model, indeed, diverges (results not shown). 

\subsection{Data-driven computation of the CLVs through the FEM-BV-VAR}

The direct computation of CLVs requires an analytical expression of the linearized dynamics (in order to apply Froyland's algorithms to the linear propagator). Hence, such a computation is not feasible via purely data-driven approaches. The idea introduced by \cite{Quinn_2021} is to use the auto-regressive linear model obtained by the FEM-BV-VAR clustering step as an underlying model to describe the dynamical system. Let us recall that the FEM-BV-VAR gives us a VAR model for each of the $K$ states

$$
\mathbf{x}_{t}=\mu^{(i)}(t)+\sum_{\tau=1}^{m} \mathbf{A}_{\tau}^{(i)}(t) \mathbf{x}_{t-\tau}+\epsilon_{t}^{(i)}
$$

From this we deduce a discrete linear dynamical system (here given for $m = 3$) :
$$
\left[\begin{array}{c}
\mathbf{x}_{t+1} \\
\mathbf{x}_{t} \\
\mathbf{x}_{t-1}
\end{array}\right]=\left[\begin{array}{ccc}
\mathbf{A}_{1}^{\left(i_{t+1}\right)} & \mathbf{A}_{2}^{\left(i_{t+1}\right)} & \mathbf{A}_{3}^{\left(i_{t+1}\right)} \\
\mathbf{I} & \mathbf{0} & \mathbf{0} \\
\mathbf{0} & \mathbf{I} & \mathbf{0}
\end{array}\right]\left[\begin{array}{c}
\mathbf{x}_{t} \\
\mathbf{x}_{t-1} \\
\mathbf{x}_{t-2}
\end{array}\right]
$$
where $i_{t+1}$ is the index of the state of the system at time $t+1$. We can therefore compute the cocycle $\mathcal{F}\left(t, t+\tau\right) = \mathcal{A}(t+\tau) \ldots \mathcal{A}(t)$, with
$$
\mathcal{A}(t)=\left[\begin{array}{ccc}
\mathbf{A}_{1}^{\left(i_{t+1}\right)} & \mathbf{A}_{2}^{\left(i_{t+1}\right)} & \mathbf{A}_{3}^{\left(i_{t+1}\right)} \\
\mathbf{I} & \mathbf{0} & \mathbf{0} \\
\mathbf{0} & \mathbf{I} & \mathbf{0}
\end{array}\right]
$$

Using the described approach, Quinn et al.~\cite{Quinn_2021} analyzed the dynamics of the North Atlantic Oscillation, using daily means of the 500 hPa geopotential height as input data. The clustering framework was used to characterise the persistent states in the atmospheric circulation, and the uncovered model was used to analyse the dynamical properties of different regimes. In particular, a finite-time dimension measure for the linear dynamical system was used to characterize the instability of each regime, thereby identifying the largest dimension to be associated with a given state of the NAO, namely the blocked state. 
They also considered the most unstable CLVs just before a transition from one state to another, to investigate which atmospheric pattern was driving the instability. The results appeared consistent with previous studies based on different methodologies.
This raised the following question: to what extent are the CLVs, computed in such a manner, significant dynamic indicators and can this method be applied to a large class of systems?
In the following, we will test thoroughly this method on systems for which many dynamical aspects are known: a fast-slow FitzHugh-Nagumo oscillator, a well-studied Von Kármán turbulent flow from a laboratory experiment, and a Lorenz 63 system. 

\section{Observations and guidelines}

The purpose of the study is to determine the conditions under which the results obtained by computing the CLVs of a data series through the FEM-BV-VAR model are reliable. The method is applied to systems for which a priori knowledge of the states and of their stability exists. In terms of dynamical structure, the examples are introduced following an increase in complexity: the method is first applied on a fast-slow FitzHugh-Nagumo oscillator with two distinct time scales, then on data extracted from a laboratory experiment of a flow whose dynamics highlight a periodic orbit and a saddle point. Finally, the  chaotic Lorenz attractor, which presents the most complex dynamics, is investigated.

Our main finding is that this procedure works well provided the studied system exhibits two properties (which are related to each other). Firstly, it should have a clear scale separation in time, that is, one should be able to distinguish a time scale gap between two (or more) phenomena in the dynamics,  as, for instance, in standard fast-slow systems. Scale separation can be estimated in several different ways, depending on the availability of data and on the existence of differential equations to describe the dynamics~\cite{rodenbeck2001dynamical,wouters2013multi,shoffner2017approaches,alberti2021small}. Secondly, 
the system needs a (near-)neutral direction along trajectories which is invariant under the linear(ized) dynamics: indeed, if the system does not have any neutral direction, the angle $\theta$ is no longer a relevant quantity to evaluate the stability of a state. This condition is frequently satisfied in physical systems, exhibiting invariant center manifolds where the hyperbolic dynamics take place; these are exactly the \emph{slow manifolds} in the fast-slow situation. For the data-driven approach to be successful, this neutral direction has to be preserved by the FEM-BV-VAR reconstructed model. This is a crucial challenge as we will see in the following.

\subsection{The case of a fast-slow FitzHugh-Nagumo oscillator} \label{sec:FHN}

\begin{figure}[hbp]
\begin{subfigure}{.49\textwidth}
  \centering
  \includegraphics[width=\linewidth]{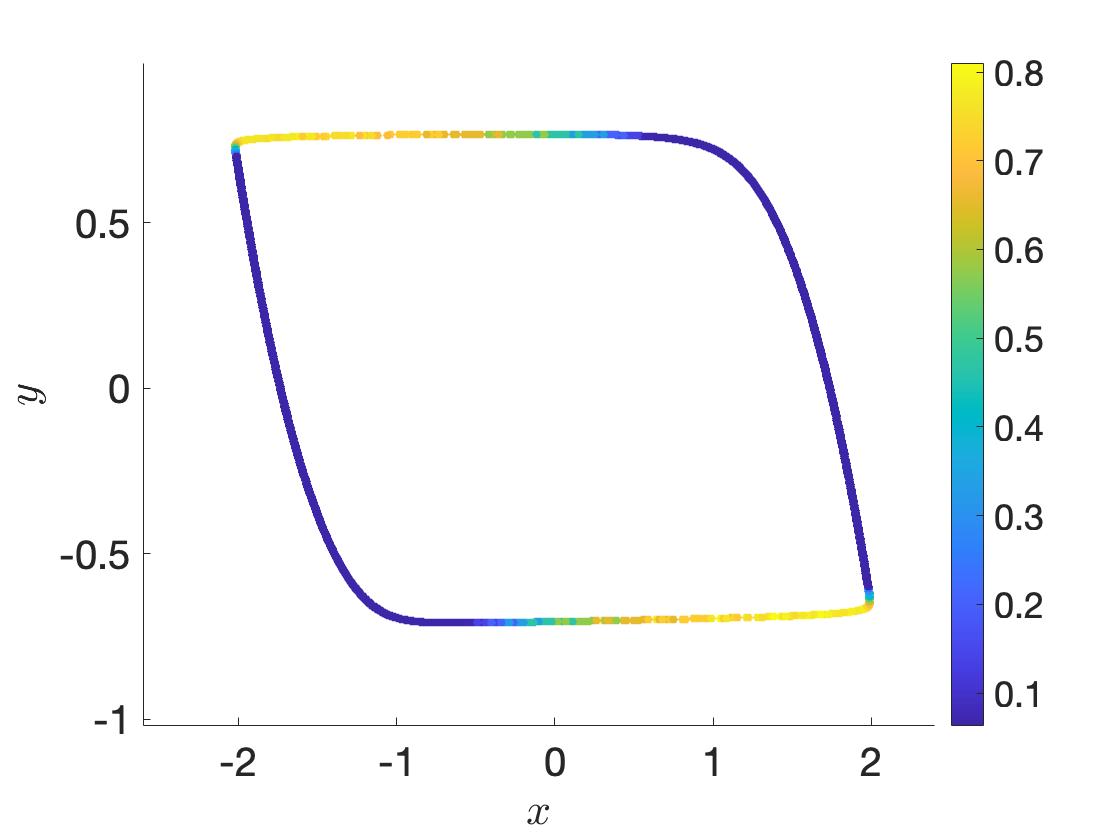}
  \caption{}
\end{subfigure}
\begin{subfigure}{.49\textwidth}
  \centering
  \includegraphics[width=\linewidth]{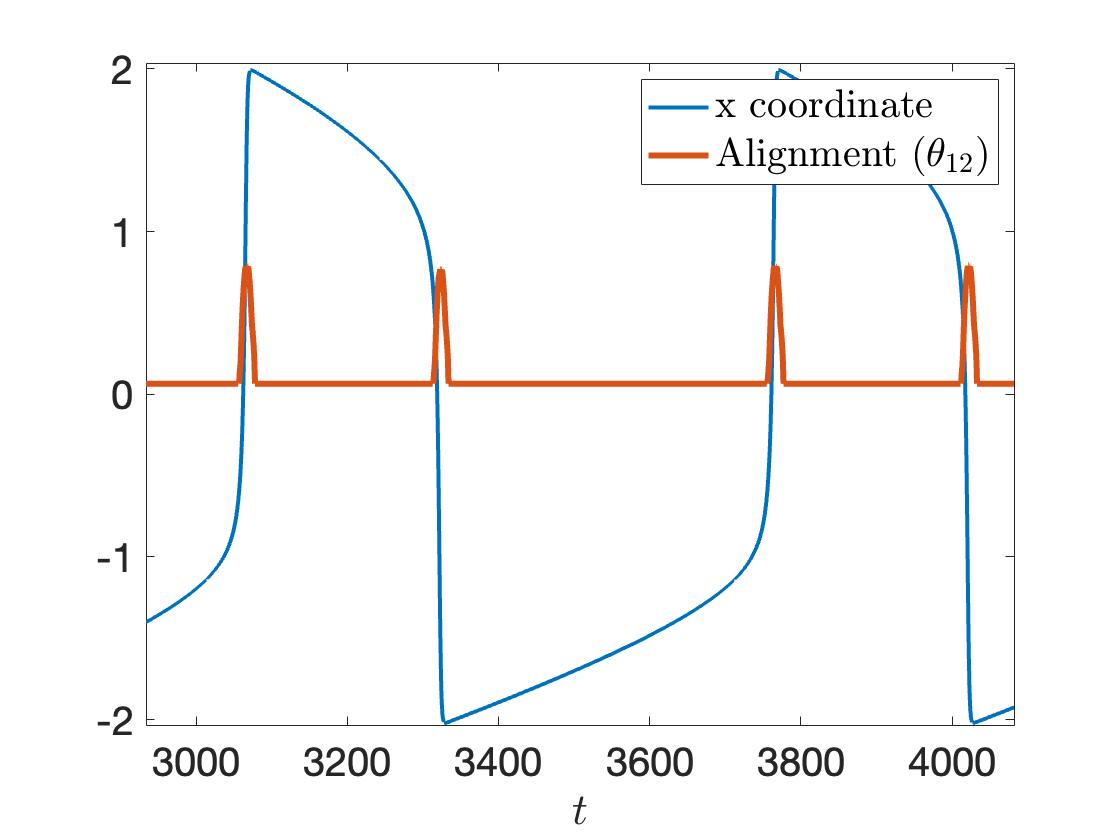}
  \caption{}
\end{subfigure}
\caption{ (a) Trajectory in the $x-y$ plane of the FitzHugh-Nagumo system, colored according to the alignment $\theta_{12}$, taking $\epsilon=0.01$, $a = 0.4$, $b = 0.3$ (standard choices, as in Sharafi et al.\cite{Sharafi_2017}). Yellow areas correspond to unstable CLVs being close to the neutral direction.  (b) Time series of the $x$ coordinate (blue) and of the alignment $\theta_{12}$ (red).}
\label{FHN}
\end{figure}

As described in Section~\ref{FEM}, the FEM-BV-VAR method is developed to study systems with a certain fast-slow structure, detecting the transition between states that are characterized by their respective fast dynamics. Therefore, the method is well-suited for models with time scale separation, expressed by a parameter $0 <\epsilon\ll 1$, that exhibit switches between different branches of the slow manifold consisting of equilibria of the fast subsystem. A by now canonical example of such a fast-slow system is the FitzHugh-Nagumo ODE~\eqref{eq:FHN} (see also Figure~\ref{FHN}), which was derived as a simplification of the Hodgkin-Huxley model for an electric potential of a nerve axon \cite{FitzHugh1955}:
\begin{align} \label{eq:FHN}
\begin{array}{r@{\;\,=\;\,}l}
\epsilon \frac{\rmd x}{\rmd \tau} & \epsilon \dot{x} = x-\frac{x^{3}}{3}-y, \\
\frac{\rmd y}{\rmd \tau} & \dot{y}  = x+a-b y.
\end{array}
\end{align}
Note that by a time change $t = \tau/\epsilon$, we may also write
\begin{align} \label{eq:FHN_fast}
\begin{array}{r@{\;\,=\;\,}l}
\frac{\rmd x}{\rmd t} & x' = x-\frac{x^{3}}{3}-y, \\
\frac{\rmd y}{\rmd t} & y'  = \epsilon (x+a-b y).
\end{array}
\end{align}
Setting $\epsilon = 0$ in equation~\eqref{eq:FHN_fast}, one can study the \emph{fast subsystem} for which $y$ is a bifurcation parameter and whose $y$-dependent set of equilibria is given by the curve $y=x-x^3/3$, also called \emph{critical manifold} $S_0$.  The cubic nonlinearity entails a bistable structure with two fold points that mark a change of stability of the fast subsystem. Considering one of the two (hyperbolically) stable branches of $S_0$, one may also take $\epsilon =0$ in equation~\eqref{eq:FHN} and observe how the \emph{slow subsystem} evolves along $S_0$. This gives a normal (or neutral) $y$-direction together with a hyperbolic $x$-direction, yielding, for $\epsilon >0$, two branches of a \emph{slow manifold} $S_{\epsilon}$ around the stable branches of $S_0$ with the same stability properties \cite{Fenichel}.
At the mentioned fold points this \emph{normal hyperbolicity} breaks down and fast switches occur between the two branches of the slow manifold (in accordance with the coloring in Figure~\ref{FHN} (a).)
The described behavior is also called \emph{relaxation-oscillation}, famously associated with the van der Pol oscillator as a paradigm model, for which the FitzHugh-Nagumo ODE is a slight generalization \cite{KuehnBook}. 
Summarizing, Figure~\ref{FHN} shows transitions between a left and a right branch of a slow manifold. Along each of these branches, there is an actual neutral direction complemented by a stable one for most of the time until both directions (almost) coincide into a locally unstable direction around the fold (or transition) points. Hence, the alignment variable $\theta_{12}$, where the stability of the CLVs is associated with the respective FTLEs, is an appropriate observable for detecting such transitions, see also Figure~\ref{FHN} (b).

In Figure~\ref{FHN}, the CLVs are computed via the FEM-BV-VAR clustering method: a FEM-BV-VAR auto-regressive model is first fitted to the timeseries of observations $(x, y)$ (see Section \ref{FEM}), for which the best hyper-parameters are found to be $K = 2$ (number of clusters), $m = 1$ (memory depth) and $p = 175$ (persistence), with an integration step $\tau = 0.003$. In this example, the choice of $K$, $m$ and $p$ is straightforward: the system has two well identifiable states, leading to $K = 2$, and the averaged persistence can easily be estimated by measuring the time spent by the system in each branch, leading to the estimate for $p$. Then, the result is fairly robust to variations in $m$, such that the simplest value $m =1$ is selected for the analysis. Having obtained an explicit linear model purely from the time series, the CLVs are approximated using the SVD-based algorithm (see Section \ref{CLVs algo}), taking $N = M = 10$ and $n = 3$. The CLV directions are robust under higher choices of $N, M$ and $n$.
Note that the sign of the associated FTLEs depends on these choices; however, since we are interested in manifesting the transition behavior happening on short time scales, the small choices of $N,M,n$ are suitable.
The alignment $\theta_{12}$ follows precisely the same profile as the one obtained through a direct computation of the CLVs from the linearization of the explicit FitzHugh-Nagumo ODE~\eqref{eq:FHN}. Sharafi et al.\cite{Sharafi_2017} also obtained a very similar pattern when they studied the CLVs of the FitzHugh-Nagumo system, based on another algorithmic procedure. Thus, the data-based method is successful for this example: via a pattern for $\theta_{12}$, one can clearly identify transitions between metastable states (corresponding with slow manifolds) through the most (finite time) unstable CLV direction (corresponding with the fast one).

The results confirm the hypothesis that systems with a clear time scale separation and a slow manifold with an actual neutral mode are well-suited for using the FEM-BV-VAR method on time series and then detecting transitions between branches of such a slow manifold  via the observable $\theta_{12}$. 

\subsection{The case of the von Kármán attractor}

Next, the method is tested on a more complex example issued from laboratory turbulent flows. In this case, the dynamics is indeed slightly more complex than in the FitzHug-Nagumo model: as will be shown in this section, an attractor can be constructed for this flow using an embedding procedure. This embedded attractor shows a periodic orbit as well as a saddle point. 

The experimental set-up is that of a von Kármán swirling flow, a device designed and maintained at the Service de Physique de l'état Condensé  of the Commissariat de l'Energie Atomique in Saclay, France~\cite{cortet2010experimental,saint2014zero,Davide_VKM,dubrulle2022many}. 
 The von Kármán turbulent flow is generated in a vertical cylinder filled with water and stirred by two coaxial, counter-rotating impellers. Those impellers provide energy and momentum flux at the upper and lower ends of the cylinder (see Fig. 2 in Dubrulle 2022\cite{dubrulle2022many}). We focus on the case where the impellers are driven by two independent motors, operating in conditions such that the torques $C_1$ and $C_2$ applied by the flow onto the top and bottom impellers are stationary. A control parameter is defined, which is capable of tracking the symmetry of the forcing, namely $\zeta = (C_1 - C_2)/(C_1+C_2)$. To quantify the global response of the flow to the forcing, the rotating frequencies $f_1$ and $f_2$ of the two impellers are measured independently. This leads to the definition of the variable $T = (f_1 - f_2)/(f_1 + f_2)$, useful to characterize the symmetries of the flow. Indeed, previous studies~\cite{saint2014zero,Davide_VKM} have identified a precise relationship between values of $T$ and instantaneous configuration of the flow: $T\simeq 0$ corresponds to a quasi-symmetric turbulent flow with two large scale circulation cells close to the impellers, and turbulence concentrated around the central section of the cylinder. For increasing $|\zeta|$, bifurcations of the flow are observed and lead to positive or negative values of $T$. Those correspond to flow geometries where a single large scale circulation structure occupies all the flow except for a turbulent boundary layer located close to the upper or lower turbine, depending on the sign of $T$. When $\abs{\zeta} > 0.06$, the von Kármán flow spontaneously switches among symmetric and bifurcated states and the dynamical switches can be approximately described by a low-dimensional attractor~\cite{Davide_VKM}. 
 
\begin{figure}[hbp]
\includegraphics[width=\linewidth]{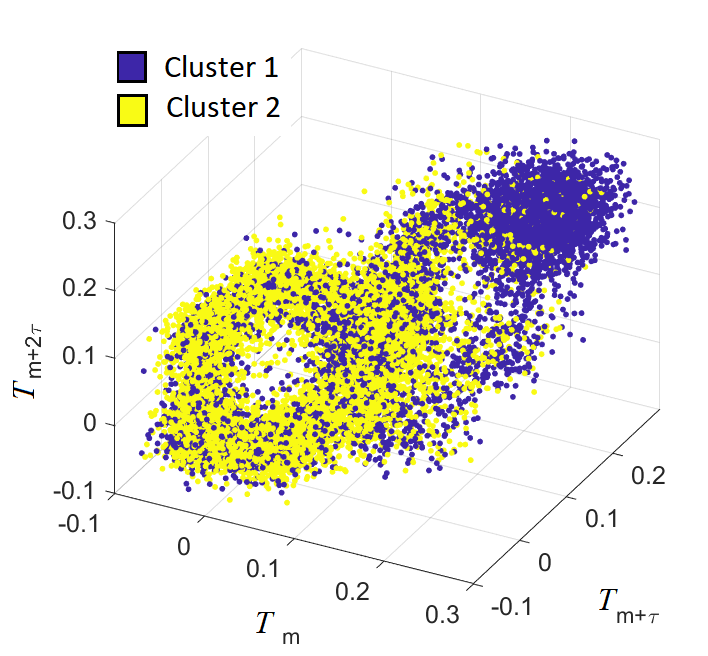}
  \caption{FEM-BV-VAR clustering on the embedded attractor for the time series of the variable $T$, from the Von Kármán experimental data. Points that the algorithm detected as part of a neighbourhood of the periodic orbit are colored in yellow, and points that are associated to the saddle point, in blue. Parameters for the FEM-BV-VAR: $K = 2$, $m = 1$, $p = 90$.}
  \label{simple VKM attractor}
\end{figure}

This attractor can be visualised with the embedding procedure, plotting $(T_m, T_{m+\tau}, T_{m+2 \tau})$. Here we will consider the case $\tau = 500$ and we refer to Faranda et al.\cite{Davide_VKM} for further details on the experiment and the choice of the parameters. The obtained embedded attractor is represented in Fig.\ref{simple VKM attractor}. It shows two persistent states: on the left a meta-stable periodic orbit, and on the right a saddle point. The system spends more time spinning around the periodic orbit than around the saddle point. From the experimental data, one can only be hypothetical about the number of unstable directions of the saddle node; however, it is clear that this fixed point supports at least one stable (attracting) and at least one unstable (repulsive) direction. We apply the FEM-BV-VAR clustering method (see Section \ref{FEM}) to the time series of $T$. To that end, the first step is to choose the best FEM-BV-VAR hyper-parameters, namely the number of states $K$, the memory depth $m$ and the persistence $p$. The embedding procedure highlights the existence of two clear states, a periodic orbit and a  saddle node, thus $K = 2$. Then a grid search is performed to select values for $m$ and $p$. As a criterion, we select the parameters that magnify the distinction between the periodic orbit and the saddle point, which corresponds to our intuition of the system behavior. The choice is based on a visual inspection of the output of the FEM-BV-VAR. The following values are finally selected: 
$$
k = 2, \quad m = 1, \quad p = 90.
$$
The corresponding state affiliation is shown in Fig.~\ref{simple VKM attractor}, where each point of the embedded attractor $(T_m, T_{m+\tau}, T_{m+2 \tau})$ is colored according to its affiliated FEM-BV-VAR cluster (also called state). One sees that the yellow state clearly corresponds to the cycle, and the blue one to the neighbourhood of the saddle point. The FEM-BV-VAR thus successfully captures the dynamical states. Let us recall that beyond the state affiliation, the FEM-BV-VAR provides a linear auto-regressive model to describe the local dynamics within a state. 

The CLVs are then computed based on the linear model given by the FEM-BV-VAR. We do not expect to have an accurate computation of the CLVs in each point, but aim at estimating the relative stability of each state.  Previous work\cite{Davide_VKM, dubrulle2022many} on the von Kármán flow experiment provide the results that can be expected: the periodic orbit is more strongly stable than the saddle point, as it is associated with the symmetric flow (see Fig. 2 in Faranda et al.\cite{Davide_VKM}). We show that the data-driven approach to compute the CLVs can retrieve this result directly from the data, looking at the alignment $\theta_{12}$ between the most unstable CLV and the near-neutral one.

\begin{figure}[htbp]
\includegraphics[width=\linewidth]{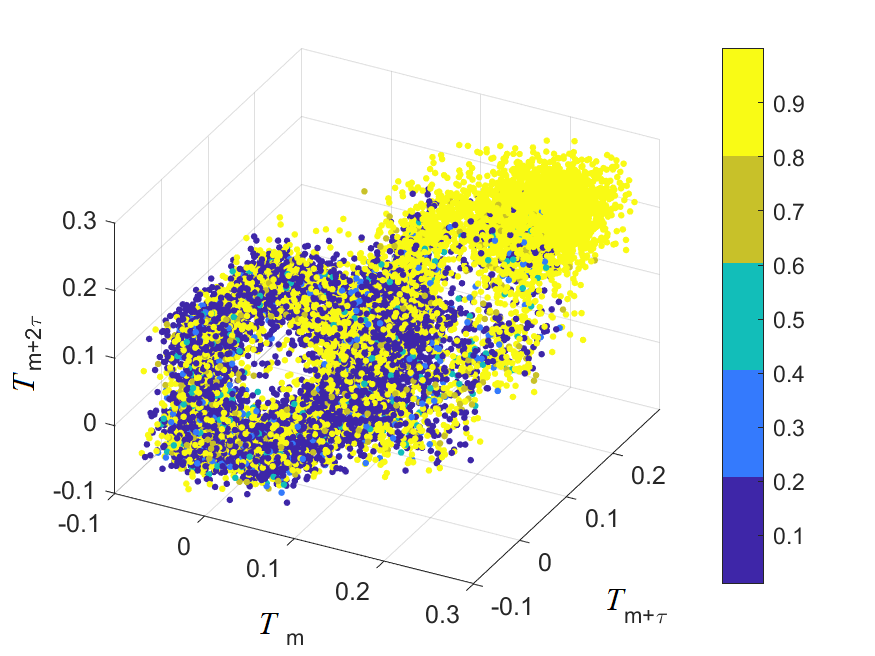}
  \caption{CLVs alignment $\theta_{12}$ on the VKM embedded attractor. Colors correspond to the value of $\theta_{12}$. In this configuration, the periodic orbit (in blue) appears to be more strongly stable than the saddle point (in yellow), which is the expected result. Parameters for the FEM-BV-VAR: $K = 2$, $m = 1$, $p = 80$. Parameters for Froyland's algorithm: $N = 30$, $n = 1$}
  \label{Align on VKM good}
\end{figure}

\begin{figure}[htbp]
\includegraphics[width=\linewidth]{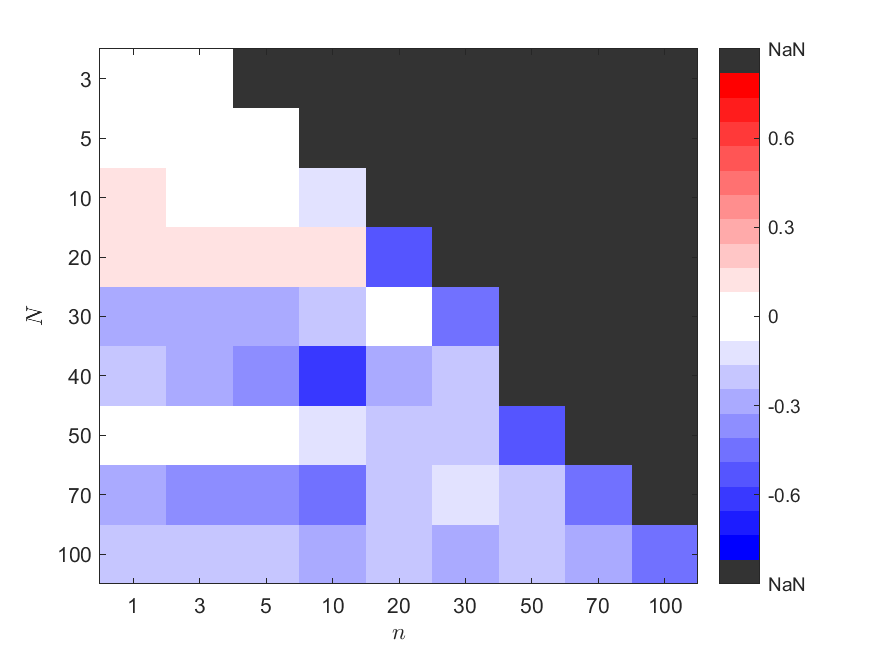}
  \caption{Difference $\Delta_{VKM}$ between the average alignment in state 1 (periodic orbit) and 2 (saddle point), as defined in Eq.~\eqref{def Delta VKM}. $N$ is the number of backward and forward steps (note that $M = N$), and $n$ is the correction step (see \ref{CLVs algo}). The blue areas correspond to the set of parameters for which the cycle is more strongly stable than the saddle point, which is expected.}
  \label{pcolor diff VKM}
\end{figure}

To that end, Froyland's algorithm (see Section \ref{CLVs algo}) is applied to the linear auto-regressive model given by the FEM-BV-VAR clustering. Three parameters need to be selected to apply the algorithm: the number of push forward steps $M$, the number of backward steps $N$ and the correction step $n$. For simplicity we take $N = M$. A grid search is then applied on $N (=M)$ and $n$. For each configuration, the CLVs and the alignment $\theta_{12}$ (as defined in Eq.~\eqref{def theta}) are computed. Fig.~\ref{Align on VKM good} shows the obtained result for one configuration of $N (=M)$ and $n$, which is consistent with the expected result. The color corresponds to the value of the alignment $\theta_{12}$, plotted on the embedded attractor, for $N = M = 30$ and $n = 1$. Around the periodic orbit the values of $\theta_{12}$ are clearly lower than around the saddle point, which means that the orbit is more strongly stable. 
However, the grid search (Fig. \ref{pcolor diff VKM}) shows that the result is not completely robust and depends on the choice of $N$ and $n$. 

To highlight the relative stability of the periodic orbit compared to the saddle point, the following difference is defined: 
\begin{equation}
    \Delta_{VKM} = \text{average of $\theta_{12}$ around the periodic orbit} - \text{average of $\theta_{12}$ around the saddle point}
    \label{def Delta VKM}
\end{equation}
Fig. \ref{pcolor diff VKM} shows, for each choice of ($N$, $n$), the value of the difference $\Delta$ between the average alignment $\theta_{12}$ on the orbit and around the saddle point. In most configurations, the difference is negative, that is to say the periodic orbit is more strongly stable than the saddle point (which is the expected result). However, care is needed because for some choices of ($N$, $n$) the result is precisely the opposite. Thus, $N$ and $n$ should be large enough, but for larger values of $N$, $\theta_{12}$ appears to become noisy (likely due to accumulation of numerical errors). Therefore the choice of $N$ and $n$ is a sensitive step, for which no systematic guidelines are available. However, the grid search used in this study supports a suitable selection of parameters, in combination with some a priori knowledge of the dynamics.

Nonetheless, this shows that for well suited values of the FEM-BV-VAR parameters (the number of states $K$, the memory depth $m$ and the persistence $p$) and of Froyland's algorithm parameters $N$ and $n$, one can obtain a very insightful information on the relative stability of the states of the system, without any \textit{a priori} information other than the raw data. This illustrates the potential validity of this method, even with experimental data. The example also supports our hypothesis that the existence of both a scale separation and a neutral direction is essential for the success of this method. In the von Kármán flow embedded attractor, one clearly has a scale separation in the sense that the trajectory oscillates for some time around one state (either the cycle or the point), and then quickly switches to the other state, with a characteristic time much faster than the oscillation. The existence of a neutral direction is more delicate to conclude, given that we do not have an underlying analytical model. However, the existence of the anticipated neutral direction is consistent with the observed quasi-periodic motion.

\subsection{On a Lorenz 63 model}
\label{sec:Lorenz}

To complete the study, the method is tested on a single Lorenz 63 system, with the usual parameters for obtaining a chaotic attractor ($\sigma = 10, \beta = 8/3, \rho = 28$)~\cite{Lorenz63}:

\begin{align}
\begin{split} \label{eq Lorenz 63}
\frac{\mathrm{d}x}{\mathrm{d}t} &= \sigma (y - x), \\
\frac{\mathrm{d}y}{\mathrm{d}t} &= x (\rho - z) - y, \\
\frac{\mathrm{d}z}{\mathrm{d}t} &= x y - \beta z.
\end{split}
\end{align}

\begin{figure}[htbp]
\centering
\includegraphics[width=\linewidth]{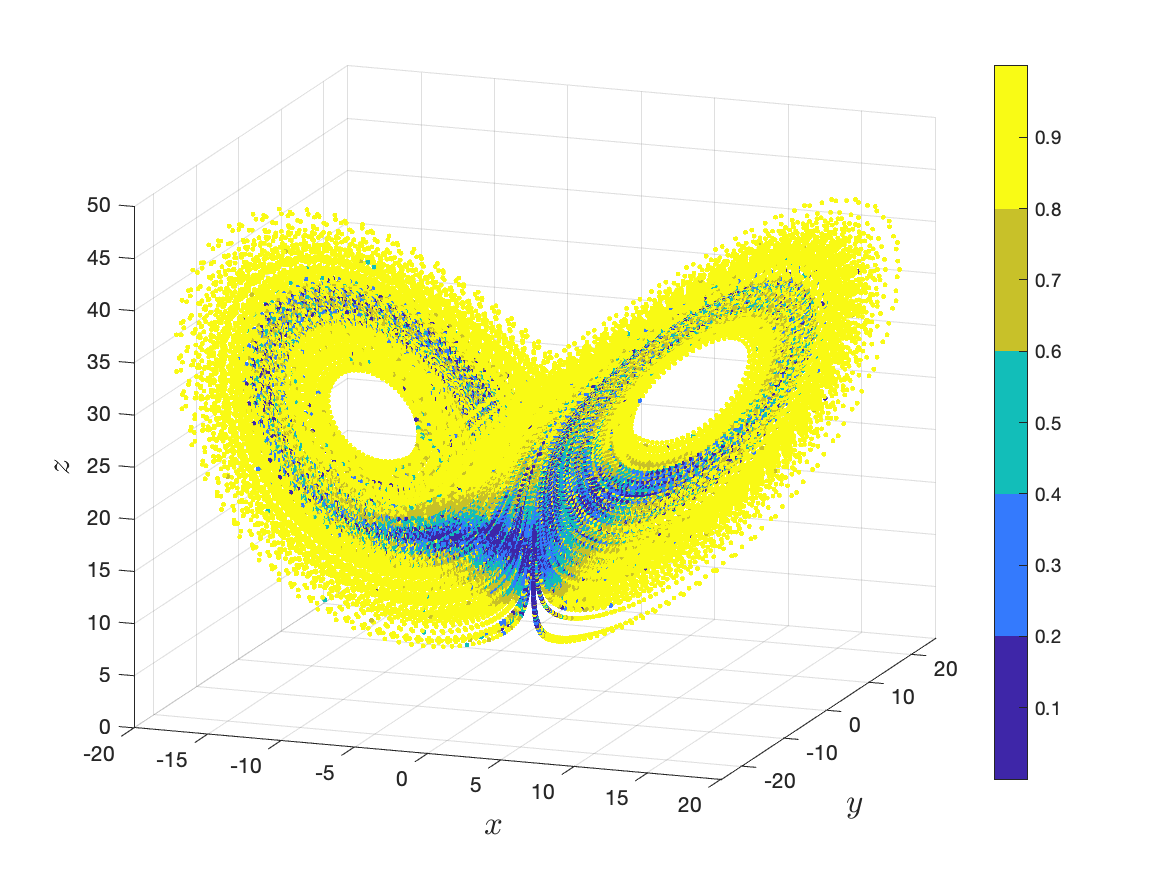}
\caption{Froyland's algorithm on a simple Lorenz 63, $N = 100$, $\tau = 0.01$ (integration step). The colors show the alignment $\theta_{12}$, as defined in Eq. \ref{def theta}. Blue areas correspond to more stable areas, where the most unstable CLV and the near-neutral one are close to being orthogonal. Conversely, yellow areas are very unstable. This result proves robust under an increase of $N$, provided $N \geq 50$.}
\label{Froyland L63 angle}
\end{figure}

\begin{figure}[hbp]
    \centering   \includegraphics[width=\textwidth]{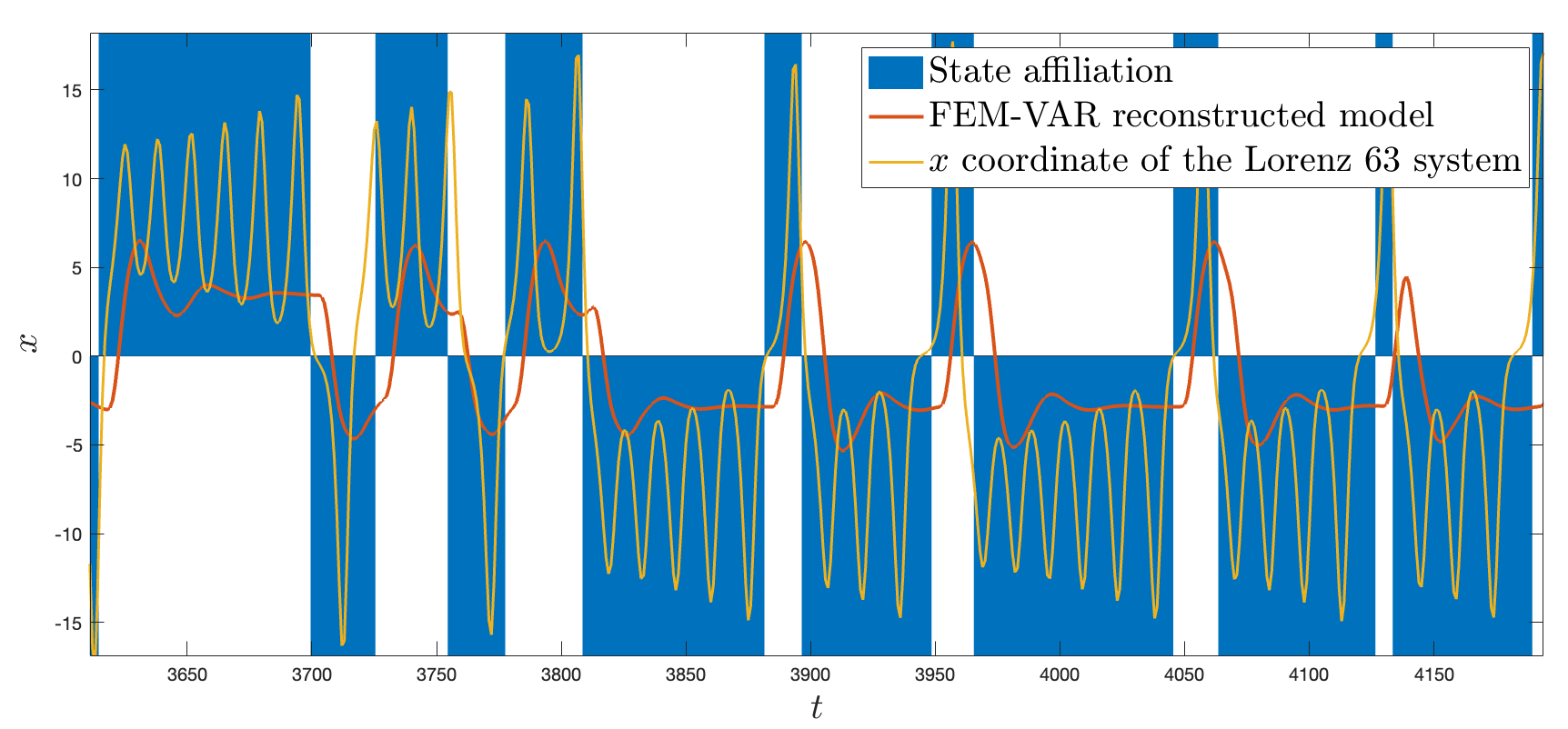}
    \caption{FEM-BV-VAR clustering applied to a Lorenz 63. First component of the Lorenz system (yellow), states affiliation (blue and white strips) and reconstruction by the FEM-BV-VAR (red). For $K = 2$, $m = 3$, $p = 29$.}
    \label{series FEM clustering L63}
\end{figure}

The attractor is self-excited with respect to three equilibria: two unstable equilibria at the center of each wing and one saddle node at the origin, see Fig. \ref{Froyland L63 angle}. The system exhibits no attracting limit cycle such that the oscillations within each wing are aperiodic, exhibiting no asymptotically exact neutral direction for the linearization.
The dynamics in each of the wings is sometimes described as metastable, with fast switches between them, such that one might think of a time scale separation. However, the associated patterns are highly irregular and not clearly associated to fast-slow dynamics (see also Figure~\ref{series FEM clustering L63}).
Dynamically speaking, this system is the most complex of this study. Regarding the Lyapunov exponents, a computation from the set of equations (\ref{eq Lorenz 63}) gives (as computed through Ginelli's procedure\cite{Ginelli_2007}):
$$
\lambda_1 = 0.9, \lambda_2 = 0.005, \lambda_3 = -14.5.
$$

These correspond to an unstable, a near-neutral and a stable direction respectively. Using the Froyland algorithm, one can compute the CLVs along the trajectory using the analytical expression of the equations (see Section \ref{CLVs algo}). Fig.~\ref{Froyland L63 angle} shows the value of the alignment $\theta_{12}$ (cosine of the angle between the most unstable CLV and the near-neutral one), plotted onto the trajectory of the Lorenz 63 system. Blue areas correspond to low values of $\theta_{12}$, therefore to more stable regions, and yellow areas to more unstable accordingly. In this study, we aim at assessing whether the FEM-BV-VAR model captures enough dynamical information for the approximated CLVs to follow a similar pattern.

As for the previous examples, one has first to choose the three parameters of the FEM-BV-VAR (namely the number of states $K$, the memory depth $m$ and the persistence $p$, see Section \ref{FEM}), which is harder in this example. $K = 2$ comes naturally as the attractor has two wings. As explained in Section \ref{choice parameters FEM}, the value of the persistence $p$ can be optimally chosen thanks to the L-curve method, provided we already fixed $K$ and $m$.  To choose $m$, the method is tested with different values of $m$ ranging from $1$ to $5$. For $m \leq 2$, the CLVs algorithm does not converge well on the FEM-BV-VAR reconstructed model. Thus we take $m = 3$, the smallest value for which the convergence is good enough. The higher $m$, the more complex the model can be (since the dimension of the auto-regressive model is $dim \times m$). With $m \leq 2$, the model may be too simple and may not capture the oscillatory patterns of the original system. Hence, the final choice is
$$
K = 2 \text{, } m = 3 \text{, } p = 29,
$$
where $p$ is chosen thanks to the L-curve method. Fig.~\ref{series FEM clustering L63} shows an extract of the time series of the original data (in yellow), the reconstructed model (in red) and the states affiliation found by the FEM-BV-VAR clustering (background in blue). 
Note that the neutral direction almost exists in the Lorenz system and leads to the oscillating dynamics. However, the FEM-BV-VAR reconstruction in Fig.~\ref{series FEM clustering L63} shows that the oscillations within a state are lost. This is a sign that the fitted AR model looses the near-neutral direction: an insight that is important for the following CLV analysis.

\begin{figure}
\includegraphics[width=\linewidth]{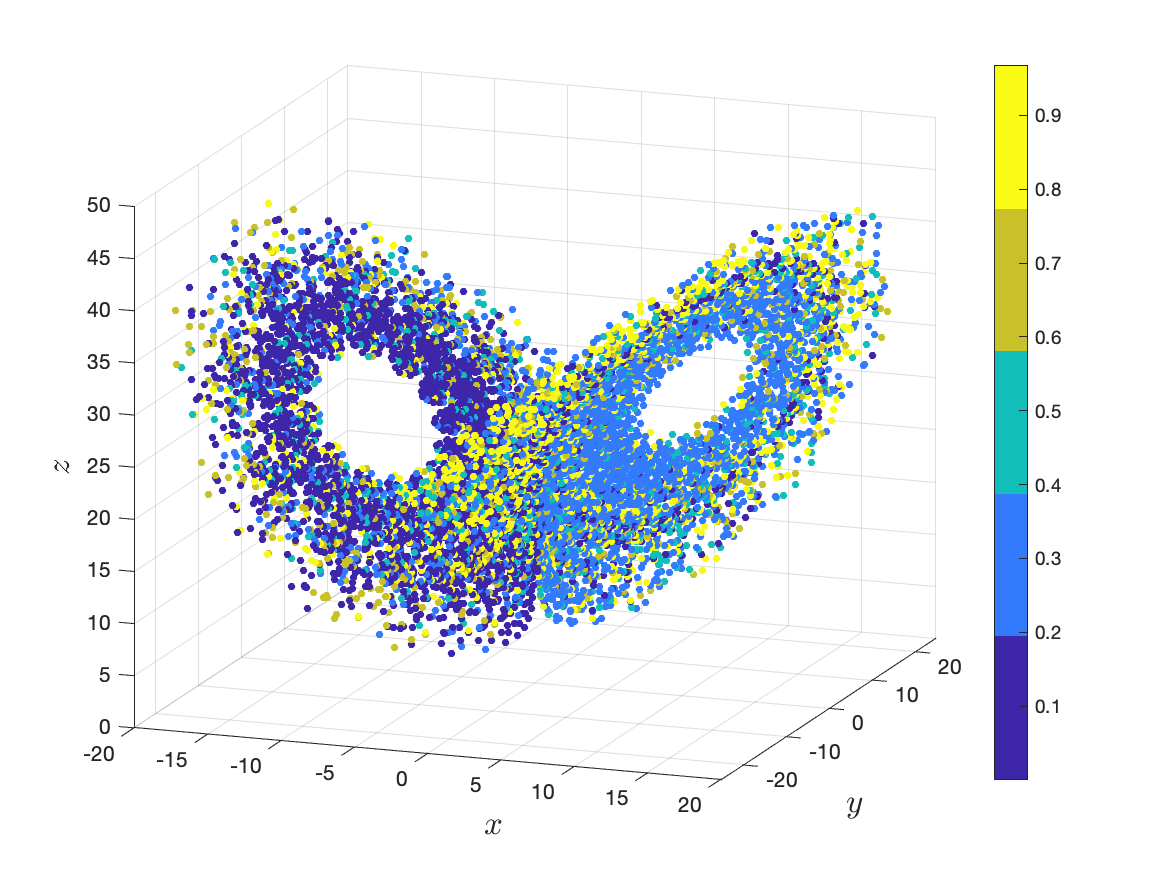}
  \caption{Alignment of CLVs on a Lorenz 63 system, obtained thanks to the FEM-BV-VAR model. In color: $\theta_{12}$. One can see that the outbound of the wings is found to be less stable than the bulk. CLVs computed with the Froyland algorithm ($N = 10$, $n = 5$), from the FEM-BV-VAR reconstruction with $K = 2$, $m = 3$, $p = 29$.}
  \label{FEM-CLVs L63 Froyland}
\end{figure}

The next step is to choose $N = M$ (the number of push backward and push forward steps) and $n$ (the correction steps) to run the CLVs algorithms (see Section \ref{CLVs algo}). It turns out that the obtained result depends highly on this choice, as for the Von Kármán flow data, except that for the Lorenz system the range of validity of the method is much narrower. For intermediate values, such as $N = 10$ and $n = 5$, one can get some information on the attractor thanks to the alignment $\theta_{12}$ obtained through the FEM-BV-VAR approach. Fig.~\ref{FEM-CLVs L63 Froyland} provides a picture that can be compared with the expected result from Fig.~\ref{Froyland L63 angle}. The absolute values of $\theta_{12}$ along the trajectories are not the same as expected. However, one can see that the outbound of the wings is found to be less stable than the bulk. Hence, the method provides again an insight on the dynamics which is, however, less precise and accurate than in the two previous examples.

To evaluate the range of validity of the method, the same picture is generated for $N$ ranging from $3$ to $100$ and $n$ from $1$ to $100$. Two criteria are used to assess the relevance of the obtained result. First, given that the distribution of the value $\theta_{12}$ has to be the same in each wing (the two wings are dynamically symmetric), the average of $\theta_{12}$ is expected to be the same in each wing. To monitor that, one can look at the difference between the average value of $\theta_{12}$ over the two wings:
\begin{equation}\label{delta lorenz}
    \Delta_{Lorenz} = \text{average of $\theta_{12}$ over the left wing} - \text{average of $\theta_{12}$ over the right wing}
\end{equation}
Secondly, to have an indicator of noisiness of the obtained time series for $\theta_{12}$, one can look at the total variation 

\begin{equation}
    TV = \sum_{i}^{} \abs{\theta_{12}(i+1) - \theta_{12}(i)}
\label{tot var}
\end{equation}

The previously shown Fig. \ref{FEM-CLVs L63 Froyland} was chosen to be the configuration that minimizes the total variation, keeping it strictly positive.

\begin{figure}[htp]
\includegraphics[width=\linewidth]{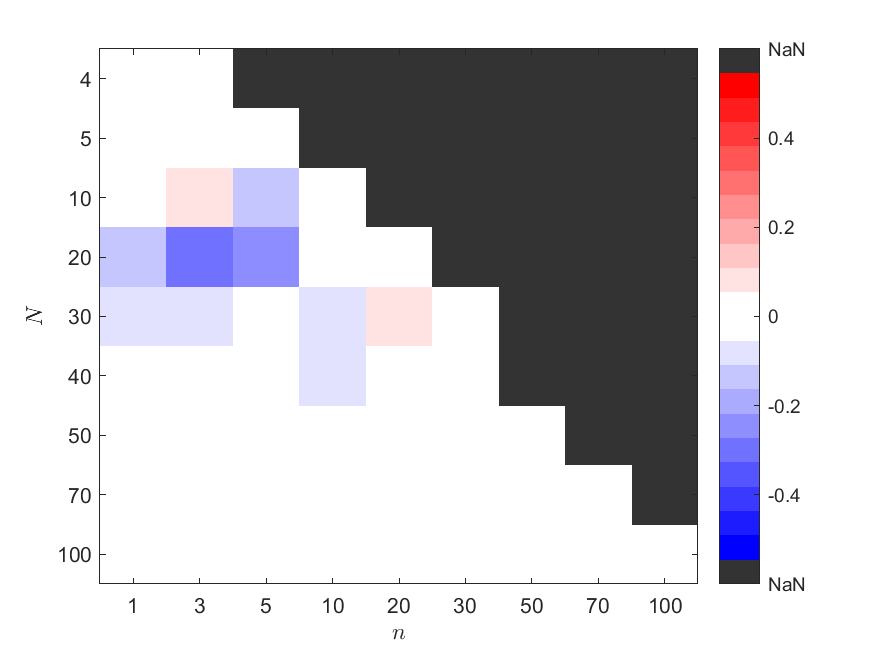}
  \caption{Pcolor plot of $\Delta_{Lorenz}$, the difference between the average value of $\theta_{12}$ over the two wings (see Eq.~\eqref{delta lorenz}). $N$ is in ordinate and $n$ in abscissa. As the two wings are symmetric, in theory this difference should be close to zero (white area). One can see that for $N$ and $n$ not large enough, $\Delta_{Lorenz}$ can be far from zero, which means that the method does not converge well with this values.}
  \label{pcolor plot sL63}
\end{figure}

Fig. \ref{pcolor plot sL63} shows, for each choice of ($N$, $n$), the value of the difference $\Delta_{Lorenz}$ between the average alignment $\theta_{12}$ over the left wing and over the right one. This value is expected to be as close as possible to zero. One can see that for small values of $N$ and $n$, the output is very asymmetric (blue zone in the bottom left), as well as for large values of $N$ (red strip on the top). As previously explained, such an asymmetry is not physically relevant. Moreover, the total variation (Eq. \eqref{tot var}) tends to increase as $N$ and $n$ increase. Thus, unlike for the von Kármán flow data, the range of validity of this method in the $N$-$n$ plane is  small, making this method hardly usable in practice for the Lorenz system. While one can have some systematic methods to tune the FEM-BV-VAR parameters ($K$, $m$ and $p$, see Section \ref{choice parameters FEM}), no such tools exist to choose $N$ and $n$. 

\begin{figure}[htp]
\begin{subfigure}{.45\textwidth}
  \centering
  \includegraphics[width=\linewidth]{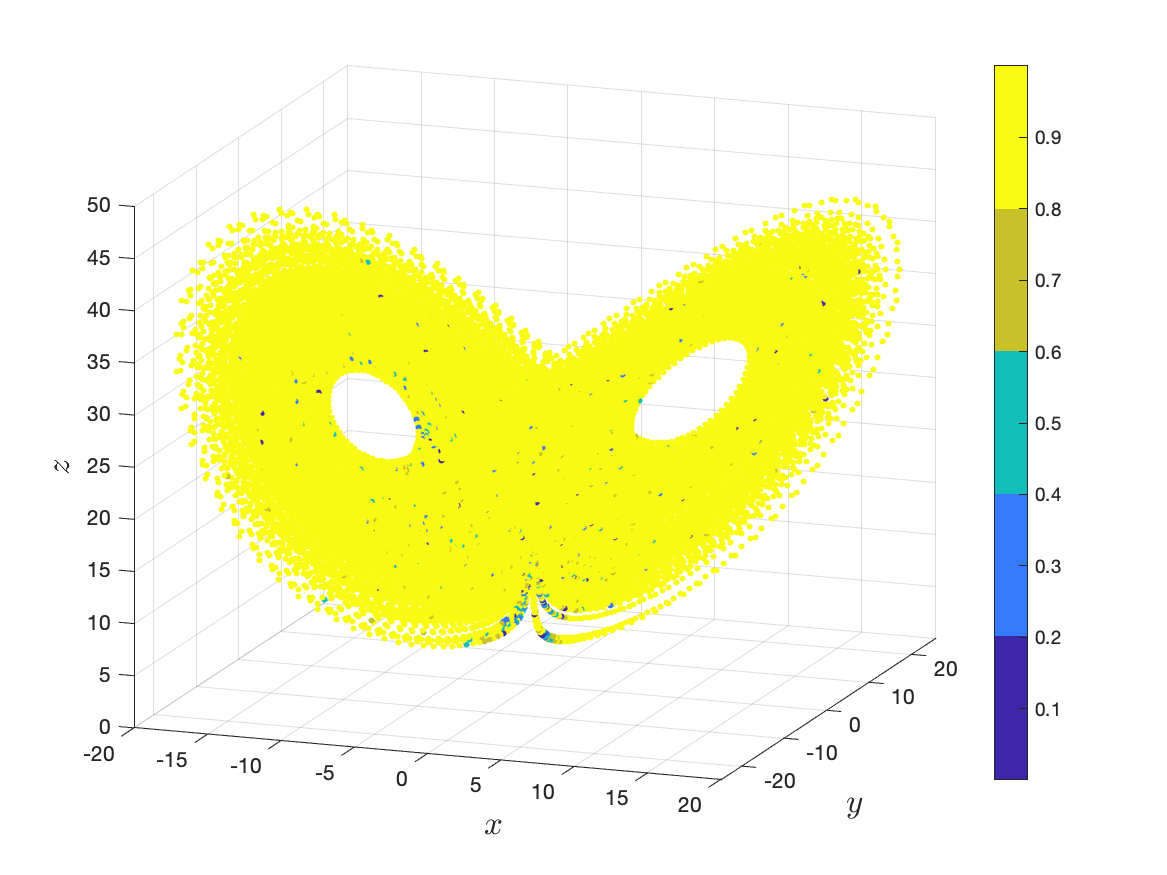}
  \caption{}
\end{subfigure}
\begin{subfigure}{.54\textwidth}
  \centering
  \includegraphics[width=\linewidth]{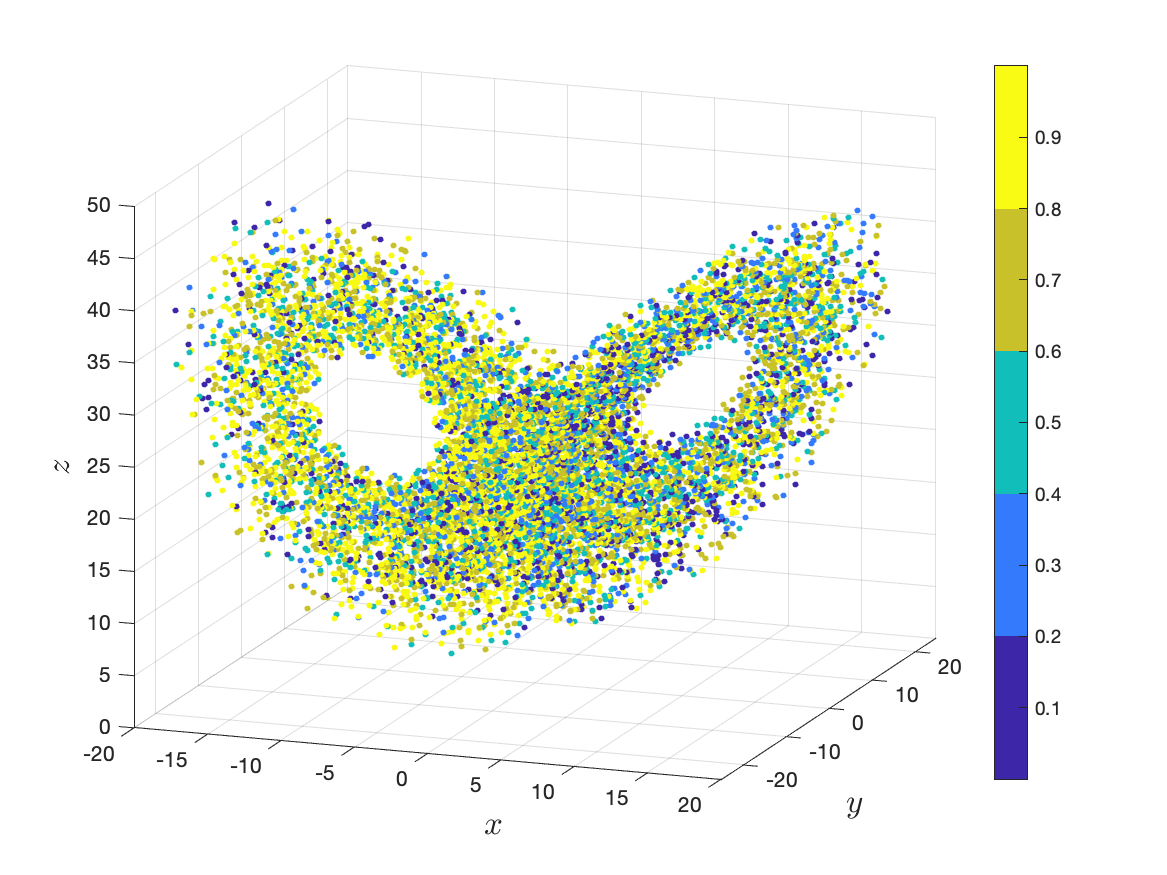}
  \caption{}
\end{subfigure}
\caption{(a) Alignment $\theta$ between the flow (tangent to the trajectory in each point) and the near neutral CLV, directly computed on a Lorenz system. Yellow corresponds to closely aligned vectors. (b) Same angle, but this time with the near neutral CLV computed on the FEM-BV-VAR reconstructed model. One can clearly see that in the first case, the near neutral CLV does correspond to the direction of the flow. However with the FEM-BV-VAR reconstruction, one completely looses this alignment.}
\label{Angle neutral flow}
\end{figure}

This observation supports our key finding: the procedure does not work well when the system has no clear time-scale separation and when the FEM-BV-VAR reconstruction does not preserve the existence of an invariant neutral direction. 
As mentioned above, one can see in Figure \ref{series FEM clustering L63} that the FEM-BV-VAR reconstructed model (in red) does not exhibit the oscillations within the wings that are characteristic of the original model (in yellow). Yet, those oscillations are important to capture the dynamics and predict the transition from one wing to the other, as suggested by Lorenz in his original paper\cite{Lorenz63}. 
In fact, the FEM-BV-VAR model seems not be able to preserve the existence of a neutral direction (of which the oscillatory dynamics are a characteristic feature). Figure \ref{Angle neutral flow} shows the alignment $\theta$ (that is to say the cosine of the angle) between the tangent to the trajectory and the expected near-neutral CLV (as there are only three dimensions, the near-neutral CLV is the second one in this case). On the left, this alignment is computed for the Lorenz 63 system directly from the analytical expression. One clearly sees that almost everywhere the second CLV and the flow are aligned, which confirms the existence of a neutral direction in this system. However, the same computation but with the CLVs computed through the FEM-BV-VAR model shows different results. The picture is completely erratic, which means that the neutral direction is (almost) entirely lost. In summary the FEM-BV-VAR model fails to capture the irregular oscillations of the system within each wing associated with such near-neutral directions. This is most likely related to the simple, linear model structure assumed in the FEM-BV-VAR approach.

\clearpage

\section{Conclusion}

The method described in this paper and suggested in earlier work by Quinn et al. \cite{Quinn_2021} makes it possible to compute an approximation of the Covariant Lyapunov Vectors (CLVs) from data series. It is based on the FEM-BV-VAR clustering scheme, which provides piece-wise auto-regressive linear models for the data. This model being built, one can compute an approximation of the CLVs. Under some conditions, the procedure seems to capture enough information on the dynamics to be able to give us the relative stability of the different areas of the phase space (that is to say, in this framework, the stability of the trajectory within each of the FEM-BV-VAR cluster). Information about stability of the trajectory is given by the analysis of the alignment between the most unstable (finite time) Lyapunov vector and the nearly neutral one (denoted $\theta_{12}$).

We claim that this procedure works well provided the studied system exhibits two properties. First, it should have a clear scale separation in time, that is to say one should be able to distinguish a temporal scale gap between two (or more) phenomena in the dynamics, as, for instance, in standard fast-slow systems. Secondly, the system has to support a dynamically invariant neutral direction in its linearization and, importantly, this neutral direction has to be preserved as much as possible by the FEM-BV-VAR reconstructed model. To support this hypothesis, we have tested the validity of the method on three different systems with an increasing dynamical complexity: the fast-slow FitzHugh-Nagumo oscillator, an embedded attractor built from von Kármán flow data that exhibits a periodic orbit along with an saddle point, and finally a classic Lorenz 63 chaotic attractor. 

In the case of the FitzHugh-Nagumo oscillator, the method yields good performances: one can find transitions precisely via the pattern of $\theta_{12}$, as the method clearly identifies switches between slow metastable regimes via unstable fast dynamics. This system exhibits a clear time-scale separation that makes it possible for the FEM-BV-VAR model to capture most of the relevant dynamical information, and especially to preserve the neutral direction. The case of the von Kármán flow shows that the method can be relevant even with experimental data, provided the dynamics exhibits a clear scale separation that allows the FEM-BV-VAR to preserve the existence of a neutral direction in the reconstructed model. It also indicates that one should be careful when tuning the values of $N$, $M$ and $n$: they must be large enough for Froyland's algorithm to converge, but not too large to avoid the accumulation of numerical errors. Finally, the Lorenz 63 example shows that for a system without a clear time scale separation, the results are highly dependent on the hyper-parameters and therefore the method is prone to fail. Due to its simple, linear model structure, the FEM-BV-VAR cannot capture irregular, complicated short term dynamical patterns (as the oscillation around the wing centers), and the reconstructed model does not show any direction that can be seen as (near-)neutral. 

Note that, while the reference approach by Quinn et al. \cite{Quinn_2021} assumes VAR models within clusters, the clustering framework introduced by Horenko \cite{Horenko_2010} is general and can accommodate more flexible model structures. Some alternative examples using different model structures can be found in Metzner et al. \cite{Metzner_2012} and in de Wiljes et al. \cite{de_wiljes_adaptive_2013}. In particular, Boyko et al. \cite{boyko_statistical_2021} recently extended this model-based clustering approach to enable the use of continuous models, effectively fitting a nonstationary, nonlinear stochastic differential equation (SDE) to timeseries of observations. Hence the data-driven computations of the CLVs could be extended to using such a SDE-based clustering for the required model fitting step. Such a future extension, based on a likely more faithful representation of complex multiscale dynamics, may lead to more accurate estimation of the CLVs and hence to a better approach to study transitions in complex systems such as the climate system.

\begin{acknowledgments}
The authors thank the Ecole Normale Superieure (ENS) for financial support enabling a research exchange of AV during which this work was started. The authors acknowledge B Dubrulle, F Daviaud and B Saint-Michel for granting the use of the von Kármán data. The work benefited from discussions with Peter Koltai. M.E. has been supported by Germany's Excellence Strategy -- The Berlin Mathematics Research Center MATH+ (EXC-2046/1, project ID: 390685689).
\end{acknowledgments}

\bibliography{Chaosarticlebib}{}

\begin{thebibliography}{37}%
\makeatletter
\providecommand \@ifxundefined [1]{%
 \@ifx{#1\undefined}
}%
\providecommand \@ifnum [1]{%
 \ifnum #1\expandafter \@firstoftwo
 \else \expandafter \@secondoftwo
 \fi
}%
\providecommand \@ifx [1]{%
 \ifx #1\expandafter \@firstoftwo
 \else \expandafter \@secondoftwo
 \fi
}%
\providecommand \natexlab [1]{#1}%
\providecommand \enquote  [1]{``#1''}%
\providecommand \bibnamefont  [1]{#1}%
\providecommand \bibfnamefont [1]{#1}%
\providecommand \citenamefont [1]{#1}%
\providecommand \href@noop [0]{\@secondoftwo}%
\providecommand \href [0]{\begingroup \@sanitize@url \@href}%
\providecommand \@href[1]{\@@startlink{#1}\@@href}%
\providecommand \@@href[1]{\endgroup#1\@@endlink}%
\providecommand \@sanitize@url [0]{\catcode `\\12\catcode `\$12\catcode
  `\&12\catcode `\#12\catcode `\^12\catcode `\_12\catcode `\%12\relax}%
\providecommand \@@startlink[1]{}%
\providecommand \@@endlink[0]{}%
\providecommand \url  [0]{\begingroup\@sanitize@url \@url }%
\providecommand \@url [1]{\endgroup\@href {#1}{\urlprefix }}%
\providecommand \urlprefix  [0]{URL }%
\providecommand \Eprint [0]{\href }%
\providecommand \doibase [0]{https://doi.org/}%
\providecommand \selectlanguage [0]{\@gobble}%
\providecommand \bibinfo  [0]{\@secondoftwo}%
\providecommand \bibfield  [0]{\@secondoftwo}%
\providecommand \translation [1]{[#1]}%
\providecommand \BibitemOpen [0]{}%
\providecommand \bibitemStop [0]{}%
\providecommand \bibitemNoStop [0]{.\EOS\space}%
\providecommand \EOS [0]{\spacefactor3000\relax}%
\providecommand \BibitemShut  [1]{\csname bibitem#1\endcsname}%
\let\auto@bib@innerbib\@empty
\bibitem [{\citenamefont {Katok}\ and\ \citenamefont
  {Hasselblatt}(1997)}]{katok1997introduction}%
  \BibitemOpen
  \bibfield  {author} {\bibinfo {author} {\bibfnamefont {A.}~\bibnamefont
  {Katok}}\ and\ \bibinfo {author} {\bibfnamefont {B.}~\bibnamefont
  {Hasselblatt}},\ }\href@noop {} {\emph {\bibinfo {title} {Introduction to the
  modern theory of dynamical systems}}},\ \bibinfo {number} {54}\ (\bibinfo
  {publisher} {Cambridge university press},\ \bibinfo {year}
  {1997})\BibitemShut {NoStop}%
\bibitem [{\citenamefont {Manneville}(2010)}]{manneville2010instabilities}%
  \BibitemOpen
  \bibfield  {author} {\bibinfo {author} {\bibfnamefont {P.}~\bibnamefont
  {Manneville}},\ }\href@noop {} {\emph {\bibinfo {title} {Instabilities, chaos
  and turbulence}}},\ Vol.~\bibinfo {volume} {1}\ (\bibinfo  {publisher} {World
  Scientific},\ \bibinfo {year} {2010})\BibitemShut {NoStop}%
\bibitem [{\citenamefont {Ruelle}(1979)}]{Ruelle_1979}%
  \BibitemOpen
  \bibfield  {author} {\bibinfo {author} {\bibfnamefont {D.}~\bibnamefont
  {Ruelle}},\ }\bibfield  {title} {\enquote {\bibinfo {title} {Sensitive
  dependence on initial condition and turbulent behavior of dynamical
  systems},}\ }\href@noop {} {\bibfield  {journal} {\bibinfo  {journal} {Annals
  of the New York Academy of Sciences}\ }\textbf {\bibinfo {volume} {316}},\
  \bibinfo {pages} {408--416} (\bibinfo {year} {1979})}\BibitemShut {NoStop}%
\bibitem [{\citenamefont {Ginelli}\ \emph {et~al.}(2007)\citenamefont
  {Ginelli}, \citenamefont {Poggi}, \citenamefont {Turchi}, \citenamefont
  {Chat{\'e}}, \citenamefont {Livi},\ and\ \citenamefont
  {Politi}}]{Ginelli_2007}%
  \BibitemOpen
  \bibfield  {author} {\bibinfo {author} {\bibfnamefont {F.}~\bibnamefont
  {Ginelli}}, \bibinfo {author} {\bibfnamefont {P.}~\bibnamefont {Poggi}},
  \bibinfo {author} {\bibfnamefont {A.}~\bibnamefont {Turchi}}, \bibinfo
  {author} {\bibfnamefont {H.}~\bibnamefont {Chat{\'e}}}, \bibinfo {author}
  {\bibfnamefont {R.}~\bibnamefont {Livi}},\ and\ \bibinfo {author}
  {\bibfnamefont {A.}~\bibnamefont {Politi}},\ }\bibfield  {title} {\enquote
  {\bibinfo {title} {Characterizing dynamics with covariant lyapunov
  vectors},}\ }\href@noop {} {\bibfield  {journal} {\bibinfo  {journal}
  {Physical review letters}\ }\textbf {\bibinfo {volume} {99}},\ \bibinfo
  {pages} {130601} (\bibinfo {year} {2007})}\BibitemShut {NoStop}%
\bibitem [{\citenamefont {Wolfe}\ and\ \citenamefont
  {Samelson}(2007)}]{Wolfe_2007}%
  \BibitemOpen
  \bibfield  {author} {\bibinfo {author} {\bibfnamefont {C.~L.}\ \bibnamefont
  {Wolfe}}\ and\ \bibinfo {author} {\bibfnamefont {R.~M.}\ \bibnamefont
  {Samelson}},\ }\bibfield  {title} {\enquote {\bibinfo {title} {An efficient
  method for recovering lyapunov vectors from singular vectors},}\ }\href@noop
  {} {\bibfield  {journal} {\bibinfo  {journal} {Tellus A: Dynamic Meteorology
  and Oceanography}\ }\textbf {\bibinfo {volume} {59}},\ \bibinfo {pages}
  {355--366} (\bibinfo {year} {2007})}\BibitemShut {NoStop}%
\bibitem [{\citenamefont {Gilmore}(2019)}]{Gilmore_2019}%
  \BibitemOpen
  \bibfield  {author} {\bibinfo {author} {\bibfnamefont {S.}~\bibnamefont
  {Gilmore}},\ }\bibfield  {title} {\enquote {\bibinfo {title} {Lyapunov
  exponents and temperature transitions in a warming australia},}\ }\href
  {https://doi.org/10.1175/JCLI-D-18-0015.1} {\bibfield  {journal} {\bibinfo
  {journal} {Journal of Climate}\ }\textbf {\bibinfo {volume} {32}},\ \bibinfo
  {pages} {2969 -- 2989} (\bibinfo {year} {2019})}\BibitemShut {NoStop}%
\bibitem [{\citenamefont {Nazarimehr}\ \emph {et~al.}(2017)\citenamefont
  {Nazarimehr}, \citenamefont {Jafari}, \citenamefont {Golpayegani},\ and\
  \citenamefont {Sprott}}]{Nazarimehr_2017}%
  \BibitemOpen
  \bibfield  {author} {\bibinfo {author} {\bibfnamefont {F.}~\bibnamefont
  {Nazarimehr}}, \bibinfo {author} {\bibfnamefont {S.}~\bibnamefont {Jafari}},
  \bibinfo {author} {\bibfnamefont {S.~M. R.~H.}\ \bibnamefont {Golpayegani}},\
  and\ \bibinfo {author} {\bibfnamefont {J.}~\bibnamefont {Sprott}},\
  }\bibfield  {title} {\enquote {\bibinfo {title} {Can lyapunov exponent
  predict critical transitions in biological systems?}}\ }\href@noop {}
  {\bibfield  {journal} {\bibinfo  {journal} {Nonlinear Dynamics}\ }\textbf
  {\bibinfo {volume} {88}},\ \bibinfo {pages} {1493--1500} (\bibinfo {year}
  {2017})}\BibitemShut {NoStop}%
\bibitem [{\citenamefont {Toth}\ and\ \citenamefont
  {Kalnay}(1993)}]{Toth_1993}%
  \BibitemOpen
  \bibfield  {author} {\bibinfo {author} {\bibfnamefont {Z.}~\bibnamefont
  {Toth}}\ and\ \bibinfo {author} {\bibfnamefont {E.}~\bibnamefont {Kalnay}},\
  }\bibfield  {title} {\enquote {\bibinfo {title} {Ensemble forecasting at nmc:
  The generation of perturbations},}\ }\href@noop {} {\bibfield  {journal}
  {\bibinfo  {journal} {Bulletin of the american meteorological society}\
  }\textbf {\bibinfo {volume} {74}},\ \bibinfo {pages} {2317--2330} (\bibinfo
  {year} {1993})}\BibitemShut {NoStop}%
\bibitem [{\citenamefont {Sharafi}, \citenamefont {Timme},\ and\ \citenamefont
  {Hallerberg}(2017)}]{Sharafi_2017}%
  \BibitemOpen
  \bibfield  {author} {\bibinfo {author} {\bibfnamefont {N.}~\bibnamefont
  {Sharafi}}, \bibinfo {author} {\bibfnamefont {M.}~\bibnamefont {Timme}},\
  and\ \bibinfo {author} {\bibfnamefont {S.}~\bibnamefont {Hallerberg}},\
  }\bibfield  {title} {\enquote {\bibinfo {title} {Critical transitions and
  perturbation growth directions},}\ }\href@noop {} {\bibfield  {journal}
  {\bibinfo  {journal} {Physical Review E}\ }\textbf {\bibinfo {volume} {96}},\
  \bibinfo {pages} {032220} (\bibinfo {year} {2017})}\BibitemShut {NoStop}%
\bibitem [{\citenamefont {Beims}\ and\ \citenamefont
  {Gallas}(2016)}]{Beims_2016}%
  \BibitemOpen
  \bibfield  {author} {\bibinfo {author} {\bibfnamefont {M.~W.}\ \bibnamefont
  {Beims}}\ and\ \bibinfo {author} {\bibfnamefont {J.~A.}\ \bibnamefont
  {Gallas}},\ }\bibfield  {title} {\enquote {\bibinfo {title} {Alignment of
  lyapunov vectors: A quantitative criterion to predict catastrophes?}}\
  }\href@noop {} {\bibfield  {journal} {\bibinfo  {journal} {Scientific
  reports}\ }\textbf {\bibinfo {volume} {6}},\ \bibinfo {pages} {1--7}
  (\bibinfo {year} {2016})}\BibitemShut {NoStop}%
\bibitem [{\citenamefont {Quinn}, \citenamefont {O'Kane},\ and\ \citenamefont
  {Kitsios}(2020)}]{Quinn_2020}%
  \BibitemOpen
  \bibfield  {author} {\bibinfo {author} {\bibfnamefont {C.}~\bibnamefont
  {Quinn}}, \bibinfo {author} {\bibfnamefont {T.~J.}\ \bibnamefont {O'Kane}},\
  and\ \bibinfo {author} {\bibfnamefont {V.}~\bibnamefont {Kitsios}},\
  }\bibfield  {title} {\enquote {\bibinfo {title} {Application of a local
  attractor dimension to reduced space strongly coupled data assimilation for
  chaotic multiscale systems},}\ }\href
  {https://doi.org/10.5194/npg-27-51-2020} {\bibfield  {journal} {\bibinfo
  {journal} {Nonlinear Processes in Geophysics}\ }\textbf {\bibinfo {volume}
  {27}},\ \bibinfo {pages} {51--74} (\bibinfo {year} {2020})}\BibitemShut
  {NoStop}%
\bibitem [{\citenamefont {Quinn}, \citenamefont {Harries},\ and\ \citenamefont
  {Kane}(2021)}]{Quinn_2021}%
  \BibitemOpen
  \bibfield  {author} {\bibinfo {author} {\bibfnamefont {C.}~\bibnamefont
  {Quinn}}, \bibinfo {author} {\bibfnamefont {D.}~\bibnamefont {Harries}},\
  and\ \bibinfo {author} {\bibfnamefont {T.~J.~O.}\ \bibnamefont {Kane}},\
  }\bibfield  {title} {\enquote {\bibinfo {title} {Dynamical analysis of a
  reduced model for the north atlantic oscillation},}\ }\href
  {https://doi.org/10.1175/JAS-D-20-0282.1} {\bibfield  {journal} {\bibinfo
  {journal} {Journal of the Atmospheric Sciences}\ }\textbf {\bibinfo {volume}
  {78}},\ \bibinfo {pages} {1647 -- 1671} (\bibinfo {year} {2021})}\BibitemShut
  {NoStop}%
\bibitem [{\citenamefont {Froyland}\ \emph {et~al.}(2013)\citenamefont
  {Froyland}, \citenamefont {H\"{u}ls}, \citenamefont {Morriss},\ and\
  \citenamefont {Watson}}]{Froyland_2013}%
  \BibitemOpen
  \bibfield  {author} {\bibinfo {author} {\bibfnamefont {G.}~\bibnamefont
  {Froyland}}, \bibinfo {author} {\bibfnamefont {T.}~\bibnamefont {H\"{u}ls}},
  \bibinfo {author} {\bibfnamefont {G.~P.}\ \bibnamefont {Morriss}},\ and\
  \bibinfo {author} {\bibfnamefont {T.~M.}\ \bibnamefont {Watson}},\ }\bibfield
   {title} {\enquote {\bibinfo {title} {Computing covariant {L}yapunov vectors,
  {O}seledets vectors, and dichotomy projectors: a comparative numerical
  study},}\ }\href {https://doi.org/10.1016/j.physd.2012.12.005} {\bibfield
  {journal} {\bibinfo  {journal} {Phys. D}\ }\textbf {\bibinfo {volume}
  {247}},\ \bibinfo {pages} {18--39} (\bibinfo {year} {2013})}\BibitemShut
  {NoStop}%
\bibitem [{\citenamefont {Horenko}(2010)}]{Horenko_2010}%
  \BibitemOpen
  \bibfield  {author} {\bibinfo {author} {\bibfnamefont {I.}~\bibnamefont
  {Horenko}},\ }\bibfield  {title} {\enquote {\bibinfo {title} {On the
  identification of nonstationary factor models and their application to
  atmospheric data analysis},}\ }\href@noop {} {\bibfield  {journal} {\bibinfo
  {journal} {Journal of the Atmospheric Sciences}\ }\textbf {\bibinfo {volume}
  {67}},\ \bibinfo {pages} {1559--1574} (\bibinfo {year} {2010})}\BibitemShut
  {NoStop}%
\bibitem [{\citenamefont {Metzner}, \citenamefont {Putzig},\ and\ \citenamefont
  {Horenko}(2012)}]{Metzner_2012}%
  \BibitemOpen
  \bibfield  {author} {\bibinfo {author} {\bibfnamefont {P.}~\bibnamefont
  {Metzner}}, \bibinfo {author} {\bibfnamefont {L.}~\bibnamefont {Putzig}},\
  and\ \bibinfo {author} {\bibfnamefont {I.}~\bibnamefont {Horenko}},\
  }\bibfield  {title} {\enquote {\bibinfo {title} {Analysis of persistent
  nonstationary time series and applications},}\ }\href@noop {} {\bibfield
  {journal} {\bibinfo  {journal} {Communications in Applied Mathematics and
  Computational Science}\ }\textbf {\bibinfo {volume} {7}},\ \bibinfo {pages}
  {175--229} (\bibinfo {year} {2012})}\BibitemShut {NoStop}%
\bibitem [{\citenamefont {Oseledec}(1968)}]{MET}%
  \BibitemOpen
  \bibfield  {author} {\bibinfo {author} {\bibfnamefont {V.~I.}\ \bibnamefont
  {Oseledec}},\ }\bibfield  {title} {\enquote {\bibinfo {title} {A
  multiplicative ergodic theorem. {C}haracteristic {L}japunov, exponents of
  dynamical systems},}\ }\href@noop {} {\bibfield  {journal} {\bibinfo
  {journal} {Trudy Moskov. Mat. Ob\v{s}\v{c}.}\ }\textbf {\bibinfo {volume}
  {19}},\ \bibinfo {pages} {179--210} (\bibinfo {year} {1968})}\BibitemShut
  {NoStop}%
\bibitem [{\citenamefont {Deremble}, \citenamefont {D’Andrea},\ and\
  \citenamefont {Ghil}(2009)}]{Deremble_2009}%
  \BibitemOpen
  \bibfield  {author} {\bibinfo {author} {\bibfnamefont {B.}~\bibnamefont
  {Deremble}}, \bibinfo {author} {\bibfnamefont {F.}~\bibnamefont
  {D’Andrea}},\ and\ \bibinfo {author} {\bibfnamefont {M.}~\bibnamefont
  {Ghil}},\ }\bibfield  {title} {\enquote {\bibinfo {title} {Fixed points,
  stable manifolds, weather regimes, and their predictability},}\ }\href@noop
  {} {\bibfield  {journal} {\bibinfo  {journal} {Chaos: An Interdisciplinary
  Journal of Nonlinear Science}\ }\textbf {\bibinfo {volume} {19}},\ \bibinfo
  {pages} {043109} (\bibinfo {year} {2009})}\BibitemShut {NoStop}%
\bibitem [{\citenamefont {Hartigan}\ and\ \citenamefont
  {Wong}(1979)}]{hartigan_algorithm_1979}%
  \BibitemOpen
  \bibfield  {author} {\bibinfo {author} {\bibfnamefont {J.~A.}\ \bibnamefont
  {Hartigan}}\ and\ \bibinfo {author} {\bibfnamefont {M.~A.}\ \bibnamefont
  {Wong}},\ }\bibfield  {title} {\enquote {\bibinfo {title} {Algorithm {AS}
  136: {A} {K}-{Means} {Clustering} {Algorithm}},}\ }\href
  {https://doi.org/10.2307/2346830} {\bibfield  {journal} {\bibinfo  {journal}
  {Journal of the Royal Statistical Society. Series C (Applied Statistics)}\
  }\textbf {\bibinfo {volume} {28}},\ \bibinfo {pages} {100--108} (\bibinfo
  {year} {1979})},\ \bibinfo {note} {publisher: [Wiley, Royal Statistical
  Society]}\BibitemShut {NoStop}%
\bibitem [{\citenamefont {Rabiner}(1989)}]{rabiner_tutorial_1989}%
  \BibitemOpen
  \bibfield  {author} {\bibinfo {author} {\bibfnamefont {L.}~\bibnamefont
  {Rabiner}},\ }\bibfield  {title} {\enquote {\bibinfo {title} {A {Tutorial} on
  {Hidden} {Markov}-{Models} and {Selected} {Applications} in {Speech}
  {Recognition}},}\ }\href
  {http://gateway.webofknowledge.com/gateway/Gateway.cgi?GWVersion=2&SrcAuth=mekentosj&SrcApp=Papers&DestLinkType=FullRecord&DestApp=WOS&KeyUT=A1989U374600002}
  {\bibfield  {journal} {\bibinfo  {journal} {Proceedings of the Ieee}\
  }\textbf {\bibinfo {volume} {77}},\ \bibinfo {pages} {257--286} (\bibinfo
  {year} {1989})}\BibitemShut {NoStop}%
\bibitem [{\citenamefont {Vercauteren}\ and\ \citenamefont
  {Klein}(2015)}]{vercauteren_clustering_2015}%
  \BibitemOpen
  \bibfield  {author} {\bibinfo {author} {\bibfnamefont {N.}~\bibnamefont
  {Vercauteren}}\ and\ \bibinfo {author} {\bibfnamefont {R.}~\bibnamefont
  {Klein}},\ }\bibfield  {title} {\enquote {\bibinfo {title} {A {Clustering}
  {Method} to {Characterize} {Intermittent} {Bursts} of {Turbulence} and
  {Interaction} with {Submesomotions} in the {Stable} {Boundary} {Layer}},}\
  }\href {https://doi.org/10.1175/JAS-D-14-0115.1} {\bibfield  {journal}
  {\bibinfo  {journal} {Journal of Atmospheric Sciences}\ }\textbf {\bibinfo
  {volume} {72}},\ \bibinfo {pages} {1504--1517} (\bibinfo {year}
  {2015})}\BibitemShut {NoStop}%
\bibitem [{\citenamefont {Boyko}\ and\ \citenamefont
  {Vercauteren}(2021)}]{boyko_multiscale_2021}%
  \BibitemOpen
  \bibfield  {author} {\bibinfo {author} {\bibfnamefont {V.}~\bibnamefont
  {Boyko}}\ and\ \bibinfo {author} {\bibfnamefont {N.}~\bibnamefont
  {Vercauteren}},\ }\bibfield  {title} {\enquote {\bibinfo {title} {Multiscale
  {Shear} {Forcing} of {Turbulence} in the {Nocturnal} {Boundary} {Layer}: {A}
  {Statistical} {Analysis}},}\ }\href
  {https://doi.org/10.1007/s10546-020-00583-0} {\bibfield  {journal} {\bibinfo
  {journal} {Boundary-Layer Meteorology}\ }\textbf {\bibinfo {volume} {179}},\
  \bibinfo {pages} {43--72} (\bibinfo {year} {2021})},\ \bibinfo {note}
  {publisher: Springer Netherlands}\BibitemShut {NoStop}%
\bibitem [{\citenamefont {O’Kane}\ \emph
  {et~al.}(2013{\natexlab{a}})\citenamefont {O’Kane}, \citenamefont {Risbey},
  \citenamefont {Franzke}, \citenamefont {Horenko},\ and\ \citenamefont
  {Monselesan}}]{FEM-OKane}%
  \BibitemOpen
  \bibfield  {author} {\bibinfo {author} {\bibfnamefont {T.~J.}\ \bibnamefont
  {O’Kane}}, \bibinfo {author} {\bibfnamefont {J.~S.}\ \bibnamefont
  {Risbey}}, \bibinfo {author} {\bibfnamefont {C.}~\bibnamefont {Franzke}},
  \bibinfo {author} {\bibfnamefont {I.}~\bibnamefont {Horenko}},\ and\ \bibinfo
  {author} {\bibfnamefont {D.~P.}\ \bibnamefont {Monselesan}},\ }\bibfield
  {title} {\enquote {\bibinfo {title} {Changes in the metastability of the
  midlatitude southern hemisphere circulation and the utility of nonstationary
  cluster analysis and split-flow blocking indices as diagnostic tools},}\
  }\href@noop {} {\bibfield  {journal} {\bibinfo  {journal} {Journal of the
  atmospheric sciences}\ }\textbf {\bibinfo {volume} {70}},\ \bibinfo {pages}
  {824--842} (\bibinfo {year} {2013}{\natexlab{a}})}\BibitemShut {NoStop}%
\bibitem [{\citenamefont {O’Kane}\ \emph
  {et~al.}(2013{\natexlab{b}})\citenamefont {O’Kane}, \citenamefont {Matear},
  \citenamefont {Chamberlain}, \citenamefont {Risbey}, \citenamefont {Sloyan},\
  and\ \citenamefont {Horenko}}]{okane_decadal_2013}%
  \BibitemOpen
  \bibfield  {author} {\bibinfo {author} {\bibfnamefont {T.~J.}\ \bibnamefont
  {O’Kane}}, \bibinfo {author} {\bibfnamefont {R.~J.}\ \bibnamefont
  {Matear}}, \bibinfo {author} {\bibfnamefont {M.~A.}\ \bibnamefont
  {Chamberlain}}, \bibinfo {author} {\bibfnamefont {J.~S.}\ \bibnamefont
  {Risbey}}, \bibinfo {author} {\bibfnamefont {B.~M.}\ \bibnamefont {Sloyan}},\
  and\ \bibinfo {author} {\bibfnamefont {I.}~\bibnamefont {Horenko}},\
  }\bibfield  {title} {\enquote {\bibinfo {title} {Decadal variability in an
  {OGCM} {Southern} {Ocean}: {Intrinsic} modes, forced modes and metastable
  states},}\ }\href {https://doi.org/10.1016/j.ocemod.2013.04.009} {\bibfield
  {journal} {\bibinfo  {journal} {Ocean Modelling}\ }\textbf {\bibinfo {volume}
  {69}},\ \bibinfo {pages} {1--21} (\bibinfo {year}
  {2013}{\natexlab{b}})}\BibitemShut {NoStop}%
\bibitem [{\citenamefont {R{\"o}denbeck}, \citenamefont {Beck},\ and\
  \citenamefont {Kantz}(2001)}]{rodenbeck2001dynamical}%
  \BibitemOpen
  \bibfield  {author} {\bibinfo {author} {\bibfnamefont {C.}~\bibnamefont
  {R{\"o}denbeck}}, \bibinfo {author} {\bibfnamefont {C.}~\bibnamefont
  {Beck}},\ and\ \bibinfo {author} {\bibfnamefont {H.}~\bibnamefont {Kantz}},\
  }\bibfield  {title} {\enquote {\bibinfo {title} {Dynamical systems with time
  scale separation: averaging, stochastic modelling, and central limit
  theorems},}\ }in\ \href@noop {} {\emph {\bibinfo {booktitle} {Stochastic
  Climate Models}}}\ (\bibinfo  {publisher} {Springer},\ \bibinfo {year}
  {2001})\ pp.\ \bibinfo {pages} {189--209}\BibitemShut {NoStop}%
\bibitem [{\citenamefont {Wouters}\ and\ \citenamefont
  {Lucarini}(2013)}]{wouters2013multi}%
  \BibitemOpen
  \bibfield  {author} {\bibinfo {author} {\bibfnamefont {J.}~\bibnamefont
  {Wouters}}\ and\ \bibinfo {author} {\bibfnamefont {V.}~\bibnamefont
  {Lucarini}},\ }\bibfield  {title} {\enquote {\bibinfo {title} {Multi-level
  dynamical systems: Connecting the ruelle response theory and the mori-zwanzig
  approach},}\ }\href@noop {} {\bibfield  {journal} {\bibinfo  {journal}
  {Journal of Statistical Physics}\ }\textbf {\bibinfo {volume} {151}},\
  \bibinfo {pages} {850--860} (\bibinfo {year} {2013})}\BibitemShut {NoStop}%
\bibitem [{\citenamefont {Shoffner}\ and\ \citenamefont
  {Schnell}(2017)}]{shoffner2017approaches}%
  \BibitemOpen
  \bibfield  {author} {\bibinfo {author} {\bibfnamefont {S.}~\bibnamefont
  {Shoffner}}\ and\ \bibinfo {author} {\bibfnamefont {S.}~\bibnamefont
  {Schnell}},\ }\bibfield  {title} {\enquote {\bibinfo {title} {Approaches for
  the estimation of timescales in nonlinear dynamical systems: Timescale
  separation in enzyme kinetics as a case study},}\ }\href@noop {} {\bibfield
  {journal} {\bibinfo  {journal} {Mathematical biosciences}\ }\textbf {\bibinfo
  {volume} {287}},\ \bibinfo {pages} {122--129} (\bibinfo {year}
  {2017})}\BibitemShut {NoStop}%
\bibitem [{\citenamefont {Alberti}\ \emph {et~al.}(2021)\citenamefont
  {Alberti}, \citenamefont {Faranda}, \citenamefont {Donner}, \citenamefont
  {Caby}, \citenamefont {Carbone}, \citenamefont {Consolini}, \citenamefont
  {Dubrulle},\ and\ \citenamefont {Vaienti}}]{alberti2021small}%
  \BibitemOpen
  \bibfield  {author} {\bibinfo {author} {\bibfnamefont {T.}~\bibnamefont
  {Alberti}}, \bibinfo {author} {\bibfnamefont {D.}~\bibnamefont {Faranda}},
  \bibinfo {author} {\bibfnamefont {R.~V.}\ \bibnamefont {Donner}}, \bibinfo
  {author} {\bibfnamefont {T.}~\bibnamefont {Caby}}, \bibinfo {author}
  {\bibfnamefont {V.}~\bibnamefont {Carbone}}, \bibinfo {author} {\bibfnamefont
  {G.}~\bibnamefont {Consolini}}, \bibinfo {author} {\bibfnamefont
  {B.}~\bibnamefont {Dubrulle}},\ and\ \bibinfo {author} {\bibfnamefont
  {S.}~\bibnamefont {Vaienti}},\ }\bibfield  {title} {\enquote {\bibinfo
  {title} {Small-scale induced large-scale transitions in solar wind magnetic
  field},}\ }\href@noop {} {\bibfield  {journal} {\bibinfo  {journal} {The
  Astrophysical journal letters}\ }\textbf {\bibinfo {volume} {914}},\ \bibinfo
  {pages} {L6} (\bibinfo {year} {2021})}\BibitemShut {NoStop}%
\bibitem [{\citenamefont {FitzHugh}(1955)}]{FitzHugh1955}%
  \BibitemOpen
  \bibfield  {author} {\bibinfo {author} {\bibfnamefont {R.}~\bibnamefont
  {FitzHugh}},\ }\bibfield  {title} {\enquote {\bibinfo {title} {Mathematical
  models of threshold phenomena in the nerve membrane},}\ }\href
  {https://doi.org/10.1007/BF02477753} {\bibfield  {journal} {\bibinfo
  {journal} {Bull. Math. Biophysics}\ }\textbf {\bibinfo {volume} {17}}
  (\bibinfo {year} {1955}),\ 10.1007/BF02477753}\BibitemShut {NoStop}%
\bibitem [{\citenamefont {Fenichel}(1979)}]{Fenichel}%
  \BibitemOpen
  \bibfield  {author} {\bibinfo {author} {\bibfnamefont {N.}~\bibnamefont
  {Fenichel}},\ }\bibfield  {title} {\enquote {\bibinfo {title} {Geometric
  singular perturbation theory for ordinary differential equations},}\ }\href
  {https://doi.org/10.1016/0022-0396(79)90152-9} {\bibfield  {journal}
  {\bibinfo  {journal} {J. Differential Equations}\ }\textbf {\bibinfo {volume}
  {31}},\ \bibinfo {pages} {53--98} (\bibinfo {year} {1979})}\BibitemShut
  {NoStop}%
\bibitem [{\citenamefont {Kuehn}(2015)}]{KuehnBook}%
  \BibitemOpen
  \bibfield  {author} {\bibinfo {author} {\bibfnamefont {C.}~\bibnamefont
  {Kuehn}},\ }\href {https://doi.org/10.1007/978-3-319-12316-5} {\emph
  {\bibinfo {title} {Multiple time scale dynamics}}},\ \bibinfo {series}
  {Applied Mathematical Sciences}, Vol.\ \bibinfo {volume} {191}\ (\bibinfo
  {publisher} {Springer, Cham},\ \bibinfo {year} {2015})\ pp.\ \bibinfo {pages}
  {xiv+814}\BibitemShut {NoStop}%
\bibitem [{\citenamefont {Cortet}\ \emph {et~al.}(2010)\citenamefont {Cortet},
  \citenamefont {Chiffaudel}, \citenamefont {Daviaud},\ and\ \citenamefont
  {Dubrulle}}]{cortet2010experimental}%
  \BibitemOpen
  \bibfield  {author} {\bibinfo {author} {\bibfnamefont {P.-P.}\ \bibnamefont
  {Cortet}}, \bibinfo {author} {\bibfnamefont {A.}~\bibnamefont {Chiffaudel}},
  \bibinfo {author} {\bibfnamefont {F.}~\bibnamefont {Daviaud}},\ and\ \bibinfo
  {author} {\bibfnamefont {B.}~\bibnamefont {Dubrulle}},\ }\bibfield  {title}
  {\enquote {\bibinfo {title} {Experimental evidence of a phase transition in a
  closed turbulent flow},}\ }\href@noop {} {\bibfield  {journal} {\bibinfo
  {journal} {Physical review letters}\ }\textbf {\bibinfo {volume} {105}},\
  \bibinfo {pages} {214501} (\bibinfo {year} {2010})}\BibitemShut {NoStop}%
\bibitem [{\citenamefont {Saint-Michel}, \citenamefont {Daviaud},\ and\
  \citenamefont {Dubrulle}(2014)}]{saint2014zero}%
  \BibitemOpen
  \bibfield  {author} {\bibinfo {author} {\bibfnamefont {B.}~\bibnamefont
  {Saint-Michel}}, \bibinfo {author} {\bibfnamefont {F.}~\bibnamefont
  {Daviaud}},\ and\ \bibinfo {author} {\bibfnamefont {B.}~\bibnamefont
  {Dubrulle}},\ }\bibfield  {title} {\enquote {\bibinfo {title} {A zero-mode
  mechanism for spontaneous symmetry breaking in a turbulent von k{\'a}rm{\'a}n
  flow},}\ }\href@noop {} {\bibfield  {journal} {\bibinfo  {journal} {New
  Journal of Physics}\ }\textbf {\bibinfo {volume} {16}},\ \bibinfo {pages}
  {013055} (\bibinfo {year} {2014})}\BibitemShut {NoStop}%
\bibitem [{\citenamefont {Faranda}\ \emph {et~al.}(2017)\citenamefont
  {Faranda}, \citenamefont {Sato}, \citenamefont {Saint-Michel}, \citenamefont
  {Wiertel}, \citenamefont {Padilla}, \citenamefont {Dubrulle},\ and\
  \citenamefont {Daviaud}}]{Davide_VKM}%
  \BibitemOpen
  \bibfield  {author} {\bibinfo {author} {\bibfnamefont {D.}~\bibnamefont
  {Faranda}}, \bibinfo {author} {\bibfnamefont {Y.}~\bibnamefont {Sato}},
  \bibinfo {author} {\bibfnamefont {B.}~\bibnamefont {Saint-Michel}}, \bibinfo
  {author} {\bibfnamefont {C.}~\bibnamefont {Wiertel}}, \bibinfo {author}
  {\bibfnamefont {V.}~\bibnamefont {Padilla}}, \bibinfo {author} {\bibfnamefont
  {B.}~\bibnamefont {Dubrulle}},\ and\ \bibinfo {author} {\bibfnamefont
  {F.}~\bibnamefont {Daviaud}},\ }\bibfield  {title} {\enquote {\bibinfo
  {title} {Stochastic chaos in a turbulent swirling flow},}\ }\href
  {https://doi.org/10.1103/PhysRevLett.119.014502} {\bibfield  {journal}
  {\bibinfo  {journal} {Phys. Rev. Lett.}\ }\textbf {\bibinfo {volume} {119}},\
  \bibinfo {pages} {014502} (\bibinfo {year} {2017})}\BibitemShut {NoStop}%
\bibitem [{\citenamefont {Dubrulle}\ \emph {et~al.}(2022)\citenamefont
  {Dubrulle}, \citenamefont {Daviaud}, \citenamefont {Faranda}, \citenamefont
  {Mari{\'e}},\ and\ \citenamefont {Saint-Michel}}]{dubrulle2022many}%
  \BibitemOpen
  \bibfield  {author} {\bibinfo {author} {\bibfnamefont {B.}~\bibnamefont
  {Dubrulle}}, \bibinfo {author} {\bibfnamefont {F.}~\bibnamefont {Daviaud}},
  \bibinfo {author} {\bibfnamefont {D.}~\bibnamefont {Faranda}}, \bibinfo
  {author} {\bibfnamefont {L.}~\bibnamefont {Mari{\'e}}},\ and\ \bibinfo
  {author} {\bibfnamefont {B.}~\bibnamefont {Saint-Michel}},\ }\bibfield
  {title} {\enquote {\bibinfo {title} {How many modes are needed to predict
  climate bifurcations? lessons from an experiment},}\ }\href@noop {}
  {\bibfield  {journal} {\bibinfo  {journal} {Nonlinear Processes in
  Geophysics}\ }\textbf {\bibinfo {volume} {29}},\ \bibinfo {pages} {17--35}
  (\bibinfo {year} {2022})}\BibitemShut {NoStop}%
\bibitem [{\citenamefont {Lorenz}(1963)}]{Lorenz63}%
  \BibitemOpen
  \bibfield  {author} {\bibinfo {author} {\bibfnamefont {E.~N.}\ \bibnamefont
  {Lorenz}},\ }\bibfield  {title} {\enquote {\bibinfo {title} {Deterministic
  nonperiodic flow},}\ }\href@noop {} {\bibfield  {journal} {\bibinfo
  {journal} {Journal of atmospheric sciences}\ }\textbf {\bibinfo {volume}
  {20}},\ \bibinfo {pages} {130--141} (\bibinfo {year} {1963})}\BibitemShut
  {NoStop}%
\bibitem [{\citenamefont {de~Wiljes}, \citenamefont {Majda},\ and\
  \citenamefont {Horenko}(2013)}]{de_wiljes_adaptive_2013}%
  \BibitemOpen
  \bibfield  {author} {\bibinfo {author} {\bibfnamefont {J.}~\bibnamefont
  {de~Wiljes}}, \bibinfo {author} {\bibfnamefont {A.~J.}\ \bibnamefont
  {Majda}},\ and\ \bibinfo {author} {\bibfnamefont {I.}~\bibnamefont
  {Horenko}},\ }\bibfield  {title} {\enquote {\bibinfo {title} {An adaptive
  {Markov} chain {Monte} {Carlo} approach to time series clustering of
  processes with regime transition behavior},}\ }\href
  {https://doi.org/10.1137/120881981} {\bibfield  {journal} {\bibinfo
  {journal} {Multiscale modeling \& simulation}\ }\textbf {\bibinfo {volume}
  {11}},\ \bibinfo {pages} {415--441} (\bibinfo {year} {2013})}\BibitemShut
  {NoStop}%
\bibitem [{\citenamefont {Boyko}, \citenamefont {Krumscheid},\ and\
  \citenamefont {Vercauteren}(2021)}]{boyko_statistical_2021}%
  \BibitemOpen
  \bibfield  {author} {\bibinfo {author} {\bibfnamefont {V.}~\bibnamefont
  {Boyko}}, \bibinfo {author} {\bibfnamefont {S.}~\bibnamefont {Krumscheid}},\
  and\ \bibinfo {author} {\bibfnamefont {N.}~\bibnamefont {Vercauteren}},\
  }\bibfield  {title} {\enquote {\bibinfo {title} {Statistical learning of
  non-linear stochastic differential equations from non-stationary time-series
  using variational clustering},}\ }\href {http://arxiv.org/abs/2102.12395}
  {\bibfield  {journal} {\bibinfo  {journal} {arXiv:2102.12395 [math]}\ }
  (\bibinfo {year} {2021})},\ \bibinfo {note} {arXiv: 2102.12395}\BibitemShut
  {NoStop}%
\end{thebibliography}%


%

\end{document}